# UTILITY MAXIMIZATION WITH A STOCHASTIC CLOCK AND AN UNBOUNDED RANDOM ENDOWMENT[1]

By Gordan Žitković

*Carnegie Mellon University*

We introduce a linear space of finitely additive measures to treat the problem of optimal expected utility from consumption under a stochastic clock and an unbounded random endowment process. In this way we establish existence and uniqueness for a large class of utility-maximization problems including the classical ones of terminal wealth or consumption, as well as the problems that depend on a random time horizon or multiple consumption instances. As an example we explicitly treat the problem of maximizing the logarithmic utility of a consumption stream, where the local time of an Ornstein–Uhlenbeck process acts as a stochastic clock.

**1. Introduction.** When we speak of expected utility, we usually have one of the following two cases in mind: expected utility of consumption on a finite interval or the expected utility of terminal wealth at some future time point. These two cases correspond to two of the historically most important problem formulations in the classical calculus of variations and optimal (stochastic) control—the *Meyer formulation* $\mathbb{E}[\int_0^T L(s, x(s))\,dt] \to \max$ and the *Lagrange formulation* $\mathbb{E}[\psi(x(T))] \to \max$, where $x(\cdot)$ denotes the controlled state function or stochastic process, and $L$ and $\psi$ correspond to the optimization criteria. These formulations owe a great deal of popularity to their analytical tractability; they fit very well into the framework of the dynamic programming principle often used to tackle optimal control problems. Even though there are a number of problem formulations in the stochastic control literature that cannot be reduced to either a Meyer or a Lagrange form [see Section 2.7, pages 85–92 of Yong and Zhou (1999), for an overview of several other classes of stochastic control models], the

Received October 2003; revised March 2004.

[1]Supported in part by NSF Grant DMS-01-39911.

expected utility theory in contemporary mathematical finance seems to lag behind in this respect. The introduction of convex duality into the treatment of utility-maximization problems by Karatzas, Lehoczky and Shreve (1987) and Karatzas, Lehoczky, Shreve and Xu (1991), as well as its further development by Kramkov and Schachermayer (1999), Cvitanić, Schachermayer and Wang (2001), Karatzas and Žitković (2003) and Hugonnier and Kramkov (2004) (to list but a small subset of the existing literature) offer hope that this lag can be overcome.

This paper aims to formulate and solve a class of utility-maximization problems of the stochastic clock type in general incomplete semimartingale markets with locally bounded stock prices and a possibly unbounded random endowment process. More specifically, our objective is to provide a mathematical framework for maximizing functionals of the form $\mathbb{E}[\int_0^T U(\omega, t, c_t) \, d\kappa_t]$, where $U$ is a time and uncertainty-dependent utility function (a utility random field), $c_t$ is the consumption density process and $\kappa_t$ is an arbitrary nondecreasing right-continuous adapted process on $[0, T]$ with $\kappa_T = 1$. Two particular choices $\kappa_t = t/T$ and $\kappa_t = \mathbf{1}_{\{t=T\}}$ correspond to the familiar Meyer and Lagrange formulations of the utility-maximization problem, but there are many other financially feasible ones. The problems of maximization of the expected utility at terminal time $T$, when $T$ is a stopping time that denotes the retirement time or a default time, form a class of examples. Another class consists of problems with the compound expected utility sampled at a sequence of stopping times. Furthermore, we could model random consumption prohibition by setting $\kappa_t = \int_0^t \mathbf{1}_{\{R_u \in C\}} \, du$ for some index process $R_t$ and a set $C \subseteq \mathbb{R}$.

The notion of a stochastic clock already was presented explicitly by Goll and Kallsen (2003) (where the phrase "stochastic clock" was introduced) and implicitly in Žitković (1999, 2002) and Karatzas and Žitković (2003). Goll and Kallsen (2003) treated the case of a logarithmic utility with no random endowment process, under additional assumptions on existence of the optimal dual process. Karatzas and Žitković (2003) established existence and uniqueness of an optimal consumption process in an incomplete semimartingale market in the presence of a bounded random endowment. Their version of the stochastic clock is, however, relatively limited—it is required to be a deterministic process with no jumps on $[0, T)$. This assumption was crucial for their treatment of the problem using convex duality and is related to the existence of a cadlag version of the optimal dual process. Related to the notion of a stochastic clock is the work by Blanchet-Scalliet, El Karoui, Jeanblanc and Martellini (2003), which deals with utility-maximization on a random horizon not necessarily given by a stopping time. Also, recent work by Bouchard and Pham (2004) treated wealth-path-dependent utility-maximization. These authors used a duality relationship between the wealth processes and a suitably chosen class of dual processes viewed as optional measures on the product space $[0,T] \times \Omega$.

In the present paper we extend the existing literature in several ways. We prove existence and describe the structure of the optimal strategy under fairly unrestrictive assumptions on the financial market and the random endowment process.

First, we allow for a general stochastic clock and a general utility that satisfies the appropriate version of the requirement of reasonable elasticity given by Kramkov and Schachermayer (1999).

Second, we allow a random endowment process that is not necessarily bounded: We require only a finite upper-hedging price for the total endowment at time $t=T$. The case of a nonbounded random endowment in the utility maximization literature was considered by Hugonnier and Kramkov (2004), but only in the case of the utility of terminal wealth and using techniques different from ours. The only restriction warranting discussion is the one we place on the jumps of the stock-price process $S$. Namely, we require $S$ to be locally bounded. The reason for this requirement [not present in Karatzas and Žitković (2003), but appearing in Hugonnier and Kramkov (2004)] is that the random endowment process is no longer assumed to be bounded and the related notion of acceptability (developed only in the locally bounded setting) has to be employed.

Finally, we present an example in which we deal completely explicitly with a utility-maximization problem in an Itô process market model with constant coefficients, where the stochastic clock is the local time at 0 of an Ornstein–Uhlenbeck process. This example illustrates how uncertainties in future consumption prohibitions introduce incompleteness into the market and describes the optimal strategy to face them.

To tackle the problem of utility maximization with the stochastic clock, we cannot depend on existing techniques. We still use the convex-duality approach, but to be able to formulate and solve the dual problem, we introduce and study the properties of two new Banach spaces: consumption densities and finitely additive measures. Also, we simplify the formulation of the standard components of the convex-duality treatment by defining the dual objective function directly as the convex conjugate of the primal objective function in the suitably coupled pair of Banach spaces. In this way, the mysterious regular parts of the finitely additive counterparts of the martingale measures used in Cvitanić, Schachermayer and Wang (2001) and Karatzas and Žitković (2003) in the definition of the dual problem appear in our treatment more naturally, in an a posteriori fashion.

The paper is organized as follows. After this Introduction, Section 2 describes the model of the financial market and poses the utility-maximization problem. In Section 3 we introduce the functional-analytic setup needed for the convex-duality treatment of our optimization problem. Section 4 introduces the convex conjugate of the utility functional and states the main result. An example that admits an explicit solution is treated in Section 5. Finally, the Appendix contains the proof of our main result.

## 2. The financial market and the optimization problem.

2.1. *The stock-price process.* We consider a financial market on a finite horizon $[0,T]$, $T \in (0,\infty)$, consisting of a $d$-dimensional locally bounded semimartingale $(S_t)_{t\in[0,T]} = (S_t^1, \ldots, S_t^d)_{t\in[0,T]}$. The process $(S_t)_{t\in[0,T]}$ is defined on a stochastic base $(\Omega, \mathcal{F}, (\mathcal{F}_t)_{t\in[0,T]}, \mathbb{P})$ that satisfies the usual conditions. For simplicity we also assume that $\mathcal{F}_0$ is $\mathbb{P}$-trivial and that $\mathcal{F} = \mathcal{F}_T$. Together with the stock-price process $(S_t)_{t\in[0,T]}$, there is a numeraire asset $S^0$ and all values are denominated in terms of $S_t^0$. This amounts to the standard assumption that $(S_t^0)_{t\in[0,T]}$ is equal to the constant process 1.

2.2. *Admissible portfolio processes.* A financial agent invests in the market according to an $(\mathcal{F}_t)_{t\in[0,T]}$-predictable $S$-integrable $d$-dimensional portfolio process $(H_t)_{t\in[0,T]}$. The stochastic integral $((H \cdot S)_t)_{t\in[0,T]}$ is called the *gains process* and represents the net gains from trade for the agent who holds a portfolio with $H_t^k$ shares of the asset $k$ at time $t$, for $k = 1, \ldots, d$.

A portfolio process $(H_t)_{t\in[0,T]}$ is called *admissible* if there exists a constant $x \in \mathbb{R}$ such that $x + (H \cdot S)_t \geq 0$ for all $t \in [0,T]$, with probability 1. Furthermore, an admissible process $(H)_{t\in[0,T]}$ is called *maximal admissible* if there exists no other admissible process $(\widetilde{H})_{t\in[0,T]}$ such that

$$(H \cdot S)_T \leq (\widetilde{H} \cdot S)_T \quad \text{a.s.} \quad \text{and} \quad \mathbb{P}[(H \cdot S)_T < (\widetilde{H} \cdot S)_T] > 0.$$

The family of all processes $(X_t^H)_{t\in[0,T]}$ of the form $X_t^H \triangleq (H \cdot S)_t$, for an admissible $H$, is denoted by $\mathcal{X}$. The class of processes $(X_t^H)_{t\in[0,T]} \in \mathcal{X}$ that corresponds to maximal admissible portfolio processes $(H)_{t\in[0,T]}$ is denoted by $\mathcal{X}_{\max}$.

We complement the widespread notion of admissibility by the less known notion of acceptability [introduced by Delbaen and Schachermayer (1997)], because admissibility is not adequate for dealing with nonbounded random endowment processes, as was shown in the context of utility maximization from terminal wealth by Hugonnier and Kramkov (2004). A portfolio process $(H)_{t\in[0,T]}$ is called *acceptable* if it admits a decomposition $H = H^+ - H^-$ with $H^+$ admissible and $H^-$ maximal admissible.

2.3. *Absence of arbitrage.* To rule out the arbitrage opportunities in our market, we state the following assumption:

ASSUMPTION 2.1. There exists a probability measure $\mathbb{Q}$ on $\mathcal{F}$, equivalent to $\mathbb{P}$, such that the process $(S_t)_{t \in [0,T]}$ is a $\mathbb{Q}$-local martingale.

The celebrated paper of Delbaen and Schachermayer (1994) showed that the condition in Assumption 2.1 is equivalent to the notion of no free lunch with vanishing risk (NFLVR)—a concept closely related to and only slightly stronger than the classical notion of absence of arbitrage. The condition NFLVR is therefore widely excepted as an operational proxy for the absence of arbitrage, and the Assumption 2.1 will be in force throughout the rest of the paper.

The set of all measures $\mathbb{Q} \sim \mathbb{P}$ as in Assumption 2.1 is denoted by $\mathcal{M}$ and we refer to the elements of $\mathcal{M}$ as the equivalent local martingale measures.

2.4. *Endowment and consumption.* Apart from being allowed to invest in the market in an admissible way, the agent (a) is continuously getting funds from an exogenous source (random endowment) and (b) is allowed to consume parts of his or her wealth as time progresses. These capital in- and out-flows are modeled by nondecreasing processes $(\mathcal{E}_t)_{t \in [0,T]}$ and $(C_t)_{t \in [0,T]}$ in $\mathcal{V}$, where $\mathcal{V}$ denotes the set of all cadlag $(\mathcal{F}_t)_{t \in [0,T]}$-optional processes vanishing at 0 whose paths are of finite variation. Here and in the rest of the paper we always identify $\mathbb{P}$-indistinguishable processes without explicit mention.

The linear space $\mathcal{V}$ can be given the structure of a vector lattice by equipping it with a partial order $\preceq$, that is compatible with its linear structure: We declare

$F^1 \preceq F^2$     if the process $(F_t^2 - F_t^1)_{t \in [0,T]}$ has nondecreasing paths.

The cone of all nondecreasing processes in $\mathcal{V}$ is the *positive cone* of the vector lattice $\mathcal{V}$ and we denote it by $\mathcal{V}_+$. Also, the *total-variation* process $(|F|_t)_{t \in [0,T]} \in \mathcal{V}_+$ is associated with each $F \in \mathcal{V}$.

The process introduced in (a) above and denoted by $(\mathcal{E}_t)_{t \in [0,T]} \in \mathcal{V}_+$ represents the *random endowment*, that is, the value $\mathcal{E}_t$ at time $t \in [0,T]$ stands for the cumulative amount of endowment received by the agent during the interval $[0,t]$. The process $(\mathcal{E}_t)_{t \in [0,T]}$ is given exogenously and we assume that the agent exerts no control over it. On the other hand, the amount and distribution of the consumption is decided by the agent, and we model the agent's consumption strategy by the consumption process $(C_t)_{t \in [0,T]} \in \mathcal{V}_+$; the value $C_t$ is the cumulative amount spent on consumption throughout the interval $[0,t]$. We find it useful in later sections to interpret the processes in $\mathcal{V}_+$ as optional random measures on the Borel sets of $[0,T]$.

2.5. *Wealth dynamics.* Starting from the initial wealth of $x \in \mathbb{R}$ (which can be negative) and the endowment process $(\mathcal{E}_t)_{t \in [0,T]}$, our agent is free to choose an acceptable portfolio process $(H_t)_{t \in [0,T]}$ and a consumption process $(C_t)_{t \in [0,T]} \in \mathcal{V}_+$. These two processes play the role of system controls. The resulting wealth process $(X_t^{(x,H,C)})_{t \in [0,T]}$ is given by the wealth dynamics equation

$$(2.1) \qquad X_t^{(x,H,C)} \triangleq x + (H \cdot S)_t - C_t + \mathcal{E}_t, \qquad t \in [0,T].$$

A consumption process $(C)_{t \in [0,T]} \in \mathcal{V}_+$ is said to be $(x, \mathcal{E})$-financeable if there exists an acceptable portfolio process $(H)_{t \in [0,T]}$ such that $X_T^{(x,H,C)} \geq 0$ a.s. The class of all $(x, \mathcal{E})$-financeable consumption processes is denoted by $\mathcal{A}(x, \mathcal{E})$ or simply by $\mathcal{A}(x)$ when there is no possibility of confusion.

REMARK 2.1. The introduction of the concept of financeability, which suppresses explicit mention of the portfolio process $(H_t)_{t \in [0,T]}$, is justified later when we specify the objective (utility) function. It depends only on the consumption, not on the particular portfolio process used to finance it, so we find it useful to formulate a static version of the optimization problem in which the portfolio process $(H_t)_{t \in [0,T]}$ does not appear at all.

REMARK 2.2. The notion of financeability imposes a weak solvency restriction on the amount of wealth the agent can consume: Even though the total wealth process $(X_t^{(x,H,C)})_{t \in [0,T]}$ is allowed to take strictly negative values before time $T$, the agent must plan consumption and investment in such a way to be able to pay all debts by the end of the planning horizon with certainty. In other words, borrowing is permitted, but only against the future endowment so that there is no chance of default. With this interpretation it makes sense to allow the initial wealth $x$ to take negative values—the initial debt might very well be covered from the future endowment. Finally, we stress that our notion of financeability differs from the one introduced in El Karoui and Jeanblanc-Picqué (1998), where no borrowing was allowed. Treatment of a consumption problem with such a stringent financeability condition seems to require a set of techniques different from ours and we leave it for future research.

2.6. *A characterization of financeable consumption processes.* In the treatment of our utility-maximization problem in the main body of this paper, the so-called budget-constraint characterization of the set $\mathcal{A}(x)$ proves to be useful. The idea is to describe the financeable consumption processes in terms of a set of linear inequalities. We provide such a characterization in the following proposition under the assumption that the random variable $\mathcal{E}_T$ (which denotes the total cumulative endowment over the horizon $[0,T]$) admits an upper-hedging price, that is, $\mathcal{U}(\mathcal{E}_T) \triangleq \sup_{\mathbb{Q} \in \mathcal{M}} \mathbb{E}^{\mathbb{Q}}[\mathcal{E}_T] < \infty$.

PROPOSITION 2.2. *Suppose that the total endowment $\mathcal{E}_T$ admits an upper-hedging price, that is, $\mathcal{U}(\mathcal{E}_T) < \infty$. Then the process $(C_t)_{t \in [0,T]} \in \mathcal{V}_+$ is $(x, \mathcal{E})$-financeable if and only if*

$$(2.2) \qquad \mathbb{E}^{\mathbb{Q}}[C_T] \leq x + \mathbb{E}^{\mathbb{Q}}[\mathcal{E}_T] \qquad \forall \, \mathbb{Q} \in \mathcal{M}.$$

PROOF. Only if: Assume first that $(C_t)_{t \in [0,T]} \in \mathcal{A}(x, \mathcal{E})$ and pick an acceptable portfolio process $(H_t)_{t \in [0,T]}$ such that the wealth process $(X_t^{(x,H,C)})_{t \in [0,T]}$ defined in (2.1) satisfies $X_T^{(x,H,C)} \geq 0$ a.s. By the definition of acceptability, there exists a decomposition $H = H_+ - H_-$ into an admissible $H_+$ and a maximal admissible $H_-$ portfolio process. Let $\mathcal{M}'$ be the set of all $\mathbb{Q} \in \mathcal{M}$ such that $((H^- \cdot S)_t)_{t \in [0,T]}$ is a $\mathbb{Q}$-uniformly integrable martingale. For any $\mathbb{Q} \in \mathcal{M}$, the process $((H^+ \cdot S)_t)_{t \in [0,T]}$ is a $\mathbb{Q}$-local martingale bounded from below and, therefore, is a $\mathbb{Q}$ supermartingale. Hence, $((H \cdot S)_t)_{t \in [0,T]}$ is a $\mathbb{Q}$ supermartingale for all $\mathbb{Q} \in \mathcal{M}'$ and

$$(2.3) \quad \begin{aligned} 0 &\leq \mathbb{E}^{\mathbb{Q}}[X_T^{(x,H,C)} | \mathcal{F}_0] = x + \mathbb{E}^{\mathbb{Q}}[(H \cdot S)_T | \mathcal{F}_0] + \mathbb{E}^{\mathbb{Q}}[\mathcal{E}_T - C_T | \mathcal{F}_0] \\ &\leq x + \mathbb{E}^{\mathbb{Q}}[\mathcal{E}_T] - \mathbb{E}^{\mathbb{Q}}[C_T] \qquad \text{for all } \mathbb{Q} \in \mathcal{M}'. \end{aligned}$$

The set $\mathcal{M}'$ of all $\mathbb{Q} \in \mathcal{M}$ such that $H^- \cdot S$ is a $\mathbb{Q}$-uniformly integrable martingale is convex and dense in $\mathcal{M}$ in the total-variation norm [see Delbaen and Schachermayer (1997), Theorem 5.2]. Therefore, the claim follows from (2.3) and the density of $\mathcal{M}'$ in $\mathcal{M}$.

IF: Let $(C_t)_{t \in [0,T]} \in \mathcal{V}_+$ be a process that satisfies $\mathbb{E}_{\mathbb{Q}}[C_T] \leq x + \mathbb{E}_{\mathbb{Q}}[\mathcal{E}_T]$ for all $\mathbb{Q} \in \mathcal{M}$. Since $\mathcal{E}_T \geq 0$ admits an upper-hedging price, there exists a constant $p > 0$ and a maximal admissible portfolio process $(H_t^{\mathcal{E}})_{t \in [0,T]}$ such that $p + (H^{\mathcal{E}} \cdot S)_T \geq \mathcal{E}_T$ a.s. [see Lemma 5.13 in Delbaen and Schachermayer (1998)]. Define the process

$$F_t \triangleq \operatorname*{ess\,sup}_{\mathbb{Q} \in \mathcal{M}} \mathbb{E}_{\mathbb{Q}}[C_T - \mathcal{E}_T + p + (H^{\mathcal{E}} \cdot S)_T | \mathcal{F}_t]$$

and note that $F_0 \leq x + p$. Then $(F_t)_{t \in [0,T]}$ is a nonnegative $\mathbb{Q}$-supermartingale for all $\mathbb{Q} \in \mathcal{M}$, permitting a cadlag modification [see Kramkov (1996), Theorem 3.2]. Thus the optional decomposition theorem [see Kramkov (1996), Theorem 2.1] asserts the existence of an admissible portfolio processes $(H_t^F)_{t \in [0,T]}$ and a finite-variation process $(G_t)_{t \in [0,T]} \in \mathcal{V}_+$ such that

$$F_t = F_0 + (H^F \cdot S)_t - G_t \qquad \text{for all } t \in [0, T] \text{ a.s.}$$

If follows that $x + p + (H^F \cdot S)_T \geq C_T - \mathcal{E}_T + p + (H^{\mathcal{E}} \cdot S)_T$, so for the acceptable portfolio process $(H_t)_{t \in [0,T]}$, defined by $H_t \triangleq H_t^F - H_t^{\mathcal{E}}$, we have $x + (H \cdot S)_T - C_T + \mathcal{E}_T \geq 0$. $\square$

2.7. *The utility functional and the primal problem.* To define the objective function of our optimization problem, we need two principal ingredients: a utility random field and the stochastic clock process.

The notion of a utility random field as defined below appeared in Žitković (1999) and Karatzas and Žitković (2003), and we use it because of its flexibility and good analytic properties—there are no continuity requirements in the temporal argument and so it is well suited for our setting.

As for the notion of a stochastic clock, it models the the agent's (either endogenously or exogenously imposed) notion of the passage of time with respect to which the consumption rate is calculated and utility is accumulated. Several examples that often appear in mathematical finance are given below. Before that let us give the formal definition of the concepts involved:

DEFINITION 2.3.

1. A *utility random field* $U: \Omega \times [0,T] \times (0,\infty) \to \mathbb{R}$ is an $\mathcal{F} \otimes \mathcal{B}[0,t] \otimes \mathcal{B}(0,\infty)$-measurable function that satisfies the following conditions.
   (a) For a fixed $(\omega, t) \in \Omega \times [0,T]$, the function $x \mapsto U(\omega, t, x)$ is a utility function, that is, a strictly concave, increasing $C^1$ function that satisfies the Inada conditions
   $$\lim_{x \to 0+} U_x(\omega, t, x) = \infty \quad \text{and} \quad \lim_{x \to \infty} U_x(\omega, t, x) = 0 \qquad \text{a.s.},$$
   where $U_x(\cdot, \cdot, \cdot)$ denotes the derivative with respect to the last argument.
   (b) There are continuous, strictly decreasing (nonrandom) functions $K_i : (0,\infty) \to (0,\infty)$, $i = 1, 2$, that satisfy
   $$\limsup_{x \to \infty} \frac{K_2(x)}{K_1(x)} < \infty$$
   and constants $G < D \in \mathbb{R}$ such that we have
   $$K_1(x) \leq U_x(\omega, t, x) \leq K_2(x)$$
   for all $(\omega, t, x) \in \Omega \times [0,T] \times (0,\infty)$ and
   $$G \leq U(\omega, t, 1) \leq D$$
   for all $(t, \omega) \in [0,T] \times \Omega$.
   (c) For every optional process $(c_t)_{t \in [0,T]}$, the process $(U(\omega, t, c_t))_{t \in [0,T]}$ is optional.
   (d) Field $U$ is reasonably elastic, that is, it satisfies $\mathrm{AE}[U] < 1$, where $\mathrm{AE}[U]$ denotes the asymptotic elasticity of the random field $U$, defined by
   $$\mathrm{AE}[U] \triangleq \limsup_{x \to \infty} \bigg( \operatornamewithlimits{ess\,sup}_{(t,\omega) \in [0,T] \times \Omega} \frac{x U_x(\omega, t, x)}{U(\omega, t, x)} \bigg).$$

2. The *stochastic clock* $(\kappa_t)_{t\in[0,T]}$ is an arbitrary process in $\mathcal{V}_+$, such that $\kappa_T = 1$, a.s.

REMARK 2.3. The requirement $\kappa_T = 1$ in Definition 2.3 is a mere normalization. We impose it to be able to work with probability measures on the product space $[0,T] \times \Omega$ (see Section 3).

We are now in the position to define the notion of a utility functional which takes consumption processes as arguments and returns their expected utility. This expected utility [as defined below in assumption (2.4)] depends only on the part of the consumption process $(C_t)_{t\in[0,T]}$ that admitts a density with respect to the stochastic measure $d\kappa$, so the choice of a consumption plan with a nontrivial component singular to $d\kappa$ is clearly suboptimal. For that reason we restrict our attention only to consumption processes $(C_t)_{t\in[0,T]}$ whose trajectories are absolutely continuous with respect to $d\kappa$, that is, only processes of the form $C_t = \int_0^t c_t \, d\kappa_t$, for a nonnegative optional process $(c_t)_{t\in[0,T]}$, which we refer to as the consumption density of the consumption process $(C_t)_{t\in[0,T]}$. For simplicity, we assume that the random endowment admits a $d\kappa$ density $(e_t)_{t\in[0,T]}$ in that $\mathcal{E}_t = \int_0^t e_u \, d\kappa_u$ for all $t \in [0,T]$, a.s. This assumption is clearly not necessary since the restrictions, which the size of the random endowment places on the choice of the consumption process, depend only on the value $\mathcal{E}_T$, as we showed in Proposition 2.2. We impose it to simplify notation by having all ingredients defined as elements of the same Banach space (see Section 3).

The utility derived from a consumption process should, therefore, be viewed as a function of the consumption density $(c_t)_{t\in[0,T]}$ and we define the *utility functional* as a function on the set of optional processes:

$$(2.4) \quad \mathbf{U}(c) \triangleq \mathbb{E}\int_0^T U(\omega, t, c_t) \, d\kappa_t \qquad \text{for an optional process } (c_t)_{t\in[0,T]}.$$

To deal with the possibility of ambiguities of the form $(+\infty) - (-\infty)$ in Definition 2.3, we adopt the following convention, which is standard in the utility-maximization literature: When the integral $\mathbb{E}\int_0^T (U(\omega,t,c_t))^- \, d\kappa_t$ of the negative part $(U(\omega,t,c_t))^-$ of the integrand from (2.4) takes the value $-\infty$, we set $\mathbf{U}(c) = -\infty$. In other words, our financial agent is not inclined toward the risks that defy classification, as far as the utility random field $U$ is concerned. Finally, we add a mild technical integrability assumption on the utility functional $U$. It is easily satisfied by all our examples and it is crucial for the simplicity of the proof of Proposition 4.1.

ASSUMPTION 2.4. For any nonnegative optional process $(c_t)_{t\in[0,T]}$ such that $\mathbf{U}(c) > -\infty$ and any constant $0 < \delta < 1$ we have $\mathbf{U}(\delta c) > -\infty$.

2.8. *Examples of utility functionals.*

EXAMPLE 2.5 (Utility random fields). 1. Let $U(x)$ be a utility function that satisfies $\limsup_{x\to\infty} \frac{xU'(x)}{U(x)} < 1$. Also, suppose there exist functions $A \colon (0,\infty) \to \mathbb{R}$ and $B \colon (0,\infty) \to (0,\infty)$ such that $U(\delta x) > A(\delta) + B(\delta) U(x)$ for all $\delta > 0$ and $x > 0$. A family of examples of such utility functions is supplied by the HARA family

$$U_\gamma(x) = \begin{cases} \dfrac{x^\gamma - 1}{\gamma}, & \gamma < 1, \gamma \neq 0, \\ \log(x), & \gamma = 0. \end{cases}$$

Then the (deterministic) utility random field

$$U(\omega, t, x) = \exp(-\beta t) U_\gamma(x)$$

conforms to Definition 2.3, and satisfies Assumption 2.4.

2. If we take a finite number $n$ of $(\mathcal{F}_t)_{t \in [0,T]}$-stopping times $\tau_1, \ldots, \tau_n$, positive constants $\beta_1, \ldots, \beta_n$ and $n$ utility functions $U^1(\cdot), \ldots, U^n(\cdot)$ as in part 1 and define

$$U(\omega, t, x) = \sum_{i=1}^{n} \exp(-\beta_i t) U^i(x) \mathbf{1}_{\{t = \tau_i(\omega)\}},$$

the random field $U$ can be easily redefined on the complement of the union of the graphs of stopping times $\tau_i$, $i = 1, \ldots, n$, to yield a utility random field satisfying Assumption 2.4.

EXAMPLE 2.6 (Stochastic clocks I). 1. Set $\kappa_t = t$, for $t \leq T = 1$. The utility functional takes the form of utility of consumption $\mathbf{U}(c) = \mathbb{E} \int_0^1 U(\omega, t, c_t) \, dt$.

2. For $\kappa_t = 0$ for $t < T$ and for $\kappa_T = 1$, we are looking at the utility of terminal wealth $\mathbb{E}[U(X_T)]$, where $U(x) = U(\omega, T, x)$. Formally, we get an expression of the form $\mathbf{U}(c) = \mathbb{E}[U(\omega, T, c_T)]$, but clearly $c_T = X_T$ in all but suboptimal cases.

3. A combination $\kappa_t = t/2$ for $t < T = 1$, and $\kappa_T = 1$, of the two cases above models the utility of consumption and terminal wealth $\mathbf{U}(c) = \mathbb{E}[\int_0^1 U(\omega, t, c_t) \, dt + U(X_T)]$.

EXAMPLE 2.7 (Stochastic clocks II). 1. Let $\tau$ be an a.s. finite $(\mathcal{F}_t)_{t \in [0,T]}$-stopping time. We can think of $\tau$ as a random horizon such as retirement time or some other market-exit time. Then the stochastic clock $\kappa_t = 0$ for $t < \tau$, and $\kappa_t = 1$ for $t \geq \tau$ models the expected utility $\mathbb{E}[U(X_\tau)]$ of the wealth at a random time $\tau$. The random endowment $\mathcal{E}_\tau$ has the interpretation of the retirement package. In the case in which the random horizon $\tau$ is unbounded, it is enough to apply a deterministic time change to fall back within the reach of our framework.

REMARK 2.4. As the anonymous referee pointed out, the case of a random horizon $\tau$ given by a mere random (as opposed to a stopping) time can be included in this framework by defining $\kappa$ as the conditional distribution of $\tau$, given the filtration $(\mathcal{F}_t)_{t\in[0,T]}$, as in Blanchet-Scalliet, El Karoui, Jeanblanc and Martellini (2003).

2. Example 2.7 can be extended to go well with the utility function from part 2 of Example 2.5. For an $n$-tuple of $(\mathcal{F}_t)_{t\in[0,T]}$-stopping times, we set

$$\kappa_t = \sum_{i=1}^{n} \frac{1}{n} \mathbf{1}_{\{t \geq \tau_i\}},$$

so that

$$\mathbf{U}(c) = \frac{1}{n} \sum_{i=1}^{n} \mathbb{E}[\exp(-\beta_i \tau_i) U^i(c_{\tau_i})].$$

3. If we set $\kappa_t = 1 - \exp(-\beta t)$ for $t < \tau$ and $\kappa_t = 1$ for $t \geq \tau$, we can add consumption to part 1 of Example 2.7,

$$\mathbf{U}(c) = \mathbb{E}\left[\int_0^\tau \exp(-\beta t) U(\omega, t, c_t)\, dt + (1 - \exp(-\beta \tau)) U(X_\tau)\right],$$

which models the utility from consumption up to and remaining wealth at random time $\tau$. The possibly inconvenient factor $(1 - \exp(-\beta\tau))$ in front of the terminal utility term can be dealt with by absorbing it into the utility random field.

EXAMPLE 2.8 (Stochastic clocks IV). 1. In this example we model the situation when the agent is allowed to withdraw consumption funds only when a certain index process $R_t$ satisfies $R_t \in C$ for some Borel set $C \subseteq \mathbb{R}$. In terms of the stochastic clock $\kappa$, we have $\kappa_t = \min(\int_0^t \mathbf{1}_{\{R_t \in C\}}\, dt, 1)$. The $R_t$ could take the role of a political indicator in an unstable economy where the individual's funds are under strict control of the government. Only in periods of political stability (i.e., when $R_t \in C$) are the withdrawal constraints relaxed to allow withdrawal of funds from the bank. It should be stressed here that the time horizon in this example is not deterministic. It is given by the stopping time

$$\inf\left\{t > 0 : \int_0^t \mathbf{1}_{\{R_u \in C\}}\, du \geq 1\right\}.$$

2. An approximation to the situation in part 1 of Example 2.8 arises when we assume that the set $C$ is of the form $(-\varepsilon, \varepsilon)$ for a constant $\varepsilon > 0$. If $\varepsilon$ is small enough, the occupation time $\int_0^t \mathbf{1}_{\{R_u \in C\}}\, du$ can be well approximated by the scaled local time $\frac{1}{2\varepsilon} l_t^R$ of the process $R_t$ at 0. Thus, we may set $\kappa_t = 1 \wedge l_t^R$. An instance of such a local-time-driven example is treated explicitly in Section 5.

2.9. *The optimization problem.* Having introduced the notion of the utility functional, we turn to the statement of our central optimization problem and we call it the *primal problem*. We describe it in terms of its value function $u:\mathbb{R} \to \mathbb{R}$ as

$$(2.5) \qquad u(x) \triangleq \sup_{c \in \mathcal{A}(x)} \mathbf{U}(c), \qquad x \in \mathbb{R},$$

where $\mathcal{A}(x)$ denotes the set of all $d\kappa$ densities of $(x, \mathcal{E})$-financeable consumption processes. Since we are working exclusively with consumption processes that admit a $d\kappa$ density, no ambiguities should arise from this slight abuse of notation. To have a nontrivial optimization problem, we impose the following standard assumption:

ASSUMPTION 2.9. There exists a constant $x > 0$ such that $u(x) < \infty$.

REMARK 2.5. 1. Assumption 2.9 is, of course, nontrivial, although quite common in the literature. In general, it has to be checked on a case-by-case basis. In the particular case when the stock-price process is an Itô process on a Brownian filtration with bounded coefficients, Assumption 2.9 is satisfied when there exist constants $M > 0$ and $\lambda < 1$ such that

$$0 \le U(t, x) \le M(1 + x^\lambda) \qquad \text{for all } (t, x) \in [0, T] \times (0, \infty).$$

For reference, see Karatzas and Shreve [(1998), Remark 3.9, page 274].

2. Part 1(b) of Definition 2.3 of a utility random field implies that $\mathbf{U}(c) \in (-\infty, \infty)$ for any constant consumption process $(c_t)_{t \in [0, T]}$, that is, a process $(c_t)_{t \in [0, T]}$ such that $c_t \equiv x$ for some constant $x > 0$. It follows that $u(x) > -\infty$ for all $x > 0$.

**3. The functional-analytic setup.** In this section we introduce several linear spaces of stochastic processes and finitely additive measures. They prove indispensable in the convex-duality treatment of the optimization problem defined in (2.5).

3.1. *Some families of finitely additive measures.* Let $\mathcal{O}$ denote the $\sigma$-algebra of optional sets relative to the filtration $(\mathcal{F}_t)_{t \in [0, T]}$. A measure $\mathbb{Q}$ defined on $\mathcal{F}_T$ and absolutely continuous to $\mathbb{P}$ induces a measure $\mathbb{Q}_\kappa$ on $\mathcal{O}$ if we set

$$(3.1) \qquad \mathbb{Q}_\kappa[A] = \mathbb{E}^{\mathbb{Q}} \int_0^T \mathbf{1}_A(t, \omega) \, d\kappa_t \qquad \text{for } A \in \mathcal{O}.$$

For notational clarity, we always identify optional stochastic processes $(c_t)_{t \in [0, T]}$ and random variables $c$ defined on the product space $[0, T] \times \Omega$ measurable with respect to the optional $\sigma$-algebra $\mathcal{O}$. Thus, the measure $\mathbb{Q}_\kappa$ can be seen as acting on an optional processes by means of integration over $[0, T] \times \Omega$ in the Lebesgue sense. In that spirit we introduce the notation

$$\langle c, \mathbb{Q} \rangle \triangleq \int_{[0,T] \times \Omega} c \, d\mathbb{Q}, \tag{3.2}$$

for a measure $\mathbb{Q}$ on the optional $\sigma$-algebra $\mathcal{O}$, and an optional process $c$ whenever the defining integral exists. A useful representation of the action $\langle c, \mathbb{Q}_\kappa \rangle$ of $\mathbb{Q}_\kappa$ on an optional process $(c_t)_{t \in [0,T]}$ is given in the following proposition.

PROPOSITION 3.1. *Let $\mathbb{Q}$ be a measure on $\mathcal{F}_T$, that is, absolutely continuous with respect to $\mathbb{P}$. For a nonnegative optional process $(c_t)_{t \in [0,T]}$ we have*

$$\langle c, \mathbb{Q}_\kappa \rangle = \mathbb{E} \int_0^T c_t Y_t^\mathbb{Q} \, d\kappa_t,$$

*where $(Y_t^\mathbb{Q})_{t \in [0,T]}$ is the cadlag version of the martingale $(\mathbb{E}[\frac{d\mathbb{Q}}{d\mathbb{P}} | \mathcal{F}_t])_{t \in [0,T]}$.*

PROOF. Define a nondecreasing cadlag process $(C_t)_{t \in [0,T]}$ by $C_t \triangleq \int_0^t c_u \, d\kappa_u$. By the integration-by-parts formula we have

$$Y_\tau^\mathbb{Q} C_\tau = \int_0^\tau Y_{t-}^\mathbb{Q} \, dC_t + \int_0^\tau C_{t-} \, dY_t^\mathbb{Q} + \sum_{0 \leq t \leq \tau} \Delta Y_t^\mathbb{Q} \Delta C_t$$
$$= \int_0^\tau Y_t^\mathbb{Q} \, dC_t + \int_0^\tau C_{t-} \, dY_t^\mathbb{Q}$$

for every stopping time $\tau \leq T$. Following Protter [(1990), Theorem III.17, page 107], the process $(\int_0^t C_{u-} \, dY_u^\mathbb{Q})_{t \in [0,T]}$ is a local martingale, so we can find an increasing sequence of stopping times $(\tau_n)_{n \in \mathbb{N}}$ that satisfy $\mathbb{P}[\tau_n < T] \to 0$ as $n \to \infty$ such that $\mathbb{E} \int_0^{\tau_n} C_{t-} \, dY_t^\mathbb{Q} = 0$ for every $n \in \mathbb{N}$. Taking expectations and letting $n \to \infty$, the monotone convergence theorem implies that

$$\langle c, \mathbb{Q}^\kappa \rangle = \mathbb{E}^\mathbb{Q}[C_T] = \mathbb{E}[Y_T^\mathbb{Q} C_T] = \lim_{n \to \infty} \mathbb{E} \int_0^{\tau_n} Y_t^\mathbb{Q} \, dC_t$$
$$= \mathbb{E} \int_0^T Y_t^\mathbb{Q} \, dC_t = \mathbb{E} \int_0^T c_t Y_t^\mathbb{Q} \, d\kappa_t. \qquad \square$$

REMARK 3.1. Note that the advantage of Proposition 3.1 over an invocation of the Radon–Nikodym theorem is the fact that the version obtained by the Radon–Nikodym derivative is merely optional and not necessarily cadlag.

We define $\mathcal{M}_\kappa \triangleq \{\mathbb{Q}_\kappa : \mathbb{Q} \in \mathcal{M}\}$. The set $\mathcal{M}_\kappa$ corresponds naturally to the set of all martingale measures in our setting, and considering measures on the product space $[0,T] \times \Omega$ instead of the measures on $\mathcal{F}_T$ is indispensable for utility maximization with a stochastic clock. Most of the existing approaches to optimal consumption start with equivalent martingale measures on $\mathcal{F}_T$ and relate them to stochastic processes on $(\mathcal{F}_t)_{t \in [0,T]}$ through some process of regularization. In our setting, the generic structure of the stochastic clock $(\kappa_t)_{t \in [0,T]}$ renders such a line of attack impossible.

However, as it turns out, $\mathcal{M}_\kappa$ is too small for duality treatment of the utility maximization problem. We need to enlarge it to contain finitely additive as well as countably additive measures. To make headway with this enlargement, we consider the set of all bounded finitely additive measures $\mathbb{Q}$ on $\mathcal{O}$, such that $\mathbb{P}_\kappa[A] = 0$ implies $\mathbb{Q}[A] = 0$, and we denote this set by $\mathbf{ba}(\mathcal{O}, \mathbb{P}_\kappa)$. It is well known that $\mathbf{ba}(\mathcal{O}, \mathbb{P}_\kappa)$, supplied with the total-variation norm, constitutes a Banach space which is isometrically isomorphic to the topological dual of $\mathbb{L}^\infty(\mathcal{O}, \mathbb{P}_\kappa)$ [see Dunford and Schwartz (1988) or Bhaskara and Bhaskara (1983)]. The action of an element $\mathbb{Q} \in \mathbf{ba}(\mathcal{O}, \mathbb{P}_\kappa)$ on $c \in \mathbb{L}^\infty(\mathcal{O}, \mathbb{P}_\kappa)$ is denoted by $\langle c, \mathbb{Q} \rangle$—a notation that naturally supplements the one introduced in (3.2).

On the Banach space $\mathbf{ba}(\mathcal{O}, \mathbb{P}_\kappa)$ there is a canonical partial ordering transferred from the pointwise order of $\mathbb{L}^\infty(\mathcal{O}, \mathbb{P}_\kappa)$, equipping it with the structure of a Banach lattice. The positive orthant of $\mathbf{ba}(\mathcal{O}, \mathbb{P}_\kappa)$ is denoted by $\mathbf{ba}(\mathcal{O}, \mathbb{P}_\kappa)_+$. An element $\mathbb{Q} \in \mathbf{ba}(\mathcal{O}, \mathbb{P}_\kappa)_+$ is said to be purely finitely additive or singular if there exists no nontrivial countably additive $\mathbb{Q}' \in \mathbf{ba}(\mathcal{O}, \mathbb{P}_\kappa)_+$ such that $\mathbb{Q}'[A] \leq \mathbb{Q}[A]$ for all $\mathcal{A} \in \mathcal{O}$. It is the content of the Yosida–Hewitt decomposition [see Yosida and Hewitt (1952)] that each $\mathbb{Q} \in \mathbf{ba}(\mathcal{O}, \mathbb{P}_\kappa)_+$ can be uniquely decomposed as $\mathbb{Q} = \mathbb{Q}^r + \mathbb{Q}^s$, with $\mathbb{Q}^r, \mathbb{Q}^s \in \mathbf{ba}(\mathcal{O}, \mathbb{P}_\kappa)_+$, where $\mathbb{Q}^r$ is a $\sigma$-additive measure and $\mathbb{Q}^s$ is purely finitely additive.

Having defined the ambient space $\mathbf{ba}(\mathcal{O}, \mathbb{P}_\kappa)$, we turn our attention to the definition of the set $\mathcal{D}_\kappa$ which serves as a building block in the advertised enlargement of the set $\mathcal{M}_\kappa$. Let $(\mathcal{M}_\kappa)^\circ$ be the polar of $\mathcal{M}_\kappa$ in $\mathbb{L}^\infty(\mathcal{O}, \mathbb{P}_\kappa)$ and let $\mathcal{D}_\kappa$ be the polar of $(\mathcal{M}_\kappa)^\circ$ (the bipolar of $\mathcal{M}_\kappa$), that is,

$$(\mathcal{M}_\kappa)^\circ \triangleq \{c \in \mathbb{L}^\infty(\mathcal{O}, \mathbb{P}_\kappa) : \langle c, \mathbb{Q} \rangle \leq 1 \text{ for all } \mathbb{Q} \in \mathcal{M}_\kappa\},$$
$$\mathcal{D}_\kappa \triangleq \{\mathbb{Q} \in \mathbf{ba}(\mathcal{O}, \mathbb{P}_\kappa) : \langle c, \mathbb{Q} \rangle \leq 1 \text{ for all } c \in (\mathcal{M}_\kappa)^\circ\}$$

and we note immediately that $\mathcal{D}_\kappa \subseteq \mathbf{ba}(\mathcal{O}, \mathbb{P}_\kappa)_+$, because $(\mathcal{M}_\kappa)^\circ$ contains the negative orthant $-\mathbb{L}^\infty_+(\mathcal{O}, \mathbb{P}_\kappa)$ of $\mathbb{L}^\infty(\mathcal{O}, \mathbb{P}_\kappa)$.

Finally, for $y > 0$ we define

$$\mathcal{M}_\kappa(y) \triangleq \{\xi \mathbb{Q} : \xi \in [0,y], \mathbb{Q} \in \mathcal{M}_\kappa\} \quad \text{and} \quad \mathcal{D}_\kappa(y) \triangleq \{y\mathbb{Q} : \mathbb{Q} \in \mathcal{D}_\kappa\}.$$

Observe that $\mathcal{M}_\kappa(y) \subseteq \mathcal{D}_\kappa(y)$ for each $y \geq 0$. Even though $\mathcal{M}_\kappa(y)$ typically is a proper subset of $\mathcal{D}_\kappa(y)$ for any $y > 0$, the following proposition shows that the difference is, in a sense, small.

PROPOSITION 3.2. *For $y > 0$, $\mathcal{M}_\kappa(y)$ is $\sigma(\mathbf{ba}(\mathcal{O}, \mathbb{P}_\kappa), \mathbb{L}^\infty(\mathcal{O}, \mathbb{P}_\kappa))$ dense in $\mathcal{D}_\kappa(y)$.*

PROOF. It is enough to provide a proof in the case $y = 1$. We start by showing that $\mathcal{D}_\kappa(1)$ is contained in the $\sigma(\mathbf{ba}(\mathcal{O}, \mathbb{P}_\kappa), \mathbb{L}^\infty(\mathcal{O}, \mathbb{P}_\kappa))$ closure $\mathrm{Cl}(\mathcal{M}_\kappa - \mathbf{ba}(\mathcal{O}, \mathbb{P}_\kappa)_+)$ of the set $\mathcal{M}_\kappa - \mathbf{ba}(\mathcal{O}, \mathbb{P}_\kappa)_+$, where

$$\mathcal{M}_\kappa - \mathbf{ba}(\mathcal{O}, \mathbb{P}_\kappa)_+ \triangleq \{\mathbb{Q} - \mathbb{Q}' : \mathbb{Q} \in \mathcal{M}_\kappa, \mathbb{Q}' \in \mathbf{ba}(\mathcal{O}, \mathbb{P}_\kappa)_+\}.$$

Suppose, to the contrary, that there exists $\mathbb{Q}^* \in \mathcal{D}_\kappa(1) \setminus \mathrm{Cl}(\mathcal{M}_\kappa - \mathbf{ba}(\mathcal{O}, \mathbb{P}_\kappa)_+)$. By the Hahn–Banach theorem there exists an element $c^* \in \mathbb{L}^\infty(\mathcal{O}, \mathbb{P}_\kappa)$, and constants $a < b$ such that $\langle c^*, \mathbb{Q}^* \rangle \geq b$ and $\langle c^*, \mathbb{Q} \rangle \leq a$ for all $\mathbb{Q} \in \mathrm{Cl}(\mathcal{M}_\kappa - \mathbf{ba}(\mathcal{O}, \mathbb{P}_\kappa)_+)$. Since $\mathcal{M}_\kappa - \mathbf{ba}(\mathcal{O}, \mathbb{P}_\kappa)_+$ contains all negative elements of $\mathbf{ba}(\mathcal{O}, \mathbb{P}_\kappa)$, we conclude that $c^* \geq 0$ $\mathbb{P}_\kappa$-a.s. and so $0 \leq a$. Furthermore, the positivity of $b$ implies that $\mathbb{P}_\kappa[c^* > 0] > 0$, since the probability measures in $\mathcal{M}_\kappa$ are equivalent to $\mathbb{P}_\kappa$. Therefore, $0 < a < b$ and the random variable $\frac{1}{a} c^*$ belongs to $(\mathcal{M}_\kappa)^\circ$. It follows that $\langle c^*, \mathbb{Q}^* \rangle \leq a$, a contradiction with the fact that $\langle c^*, \mathbb{Q}^* \rangle \geq b$.

To finalize the proof we pick $\mathbb{Q} \in \mathcal{D}'_\kappa(1) \triangleq \{\mathbb{Q} \in \mathcal{D}_\kappa(1) : \langle 1, \mathbb{Q} \rangle = 1\}$ and take a directed set $A$ and a net $(\widetilde{\mathbb{Q}}_\alpha)_{\alpha \in A}$ in $\mathcal{M}_\kappa - \mathbf{ba}(\mathcal{O}, \mathbb{P}_\kappa)_+$ such that $\widetilde{\mathbb{Q}}_\alpha \to \mathbb{Q}$. Such a net exists thanks to the result of the first part of this proof. Each $\widetilde{\mathbb{Q}}_\alpha$ can be written as $\widetilde{\mathbb{Q}}_\alpha = \mathbb{Q}_\alpha^{\mathcal{M}_\kappa} - \mathbb{Q}_\alpha^+$ with $\mathbb{Q}_\alpha^{\mathcal{M}_\kappa} \in \mathcal{M}_\kappa$ and $\mathbb{Q}_\alpha^+ \in \mathbf{ba}(\mathcal{O}, \mathbb{P}_\kappa)_+$ for all $\alpha \in A$. Weak-* convergence of the net $\widetilde{\mathbb{Q}}_\alpha$ implies that $\langle 1, \mathbb{Q}_\alpha^+ \rangle \to 0$ and therefore $\mathbb{Q}_\alpha^+ \to 0$ in the norm and weak-* topologies. Thus $\mathbb{Q}^{\mathcal{M}_\kappa} \to \mathbb{Q}$ and we conclude that $\mathcal{M}_\kappa$ is dense in $\mathcal{D}'_\kappa(1)$. It follows immediately that $\mathcal{M}_\kappa(1)$ is dense in $\mathcal{D}_\kappa(1)$. □

3.2. *The space $\mathcal{V}_\kappa^\mathcal{M}$.* Let $\mathcal{V}_\kappa^\mathcal{M}$ stand for the vector space of all optional random processes $(c_t)_{t \in [0,T]}$ verifying

$$\|c\|_\mathcal{M} < \infty \quad \text{where } \|c\|_\mathcal{M} \triangleq \sup_{\mathbb{Q} \in \mathcal{M}_\kappa} \langle |c|, \mathbb{Q} \rangle.$$

It is quite clear that $\|\cdot\|_\mathcal{M}$ defines a norm on $\mathcal{V}_\kappa^\mathcal{M}$. We establish completeness in the following proposition.

PROPOSITION 3.3. *$(\mathcal{V}_\kappa^\mathcal{M}, \|\cdot\|_\mathcal{M})$ is a Banach space.*

PROOF. To prove that $\mathcal{V}_\kappa^\mathcal{M}$ is complete under $\|\cdot\|_\mathcal{M}$, we take a sequence $(c_n)_{n \in \mathbb{N}}$ in $\mathcal{V}_\kappa^\mathcal{M}$ such that $\sum_n \|c_n\|_\mathcal{M} < \infty$. Given a fixed but arbitrary $\widetilde{\mathbb{Q}}_\kappa \in \mathcal{M}_\kappa$, the inequality $\|c\|_\mathcal{M} \geq \langle |c|, \widetilde{\mathbb{Q}}_\kappa \rangle$ holds for every $c \in \mathcal{V}_\kappa^\mathcal{M}$ and thus the series $\sum_{n=1}^\infty |c_n|$ converges in $\mathbb{L}^1(\mathcal{O}, \widetilde{\mathbb{Q}}_\kappa)$. We can, therefore, find an optional process $c_0 \in \mathbb{L}^1(\widetilde{\mathbb{Q}}_\kappa, \mathcal{O})$ such that $c_0 = \lim_{n \to \infty} \sum_{k=1}^n c_k$ in $\mathbb{L}^1(\widetilde{\mathbb{Q}}_\kappa, \mathcal{O})$ and $\widetilde{\mathbb{Q}}_\kappa$-a.s.

For an arbitrary $\mathbb{Q}_\kappa \in \mathcal{M}_\kappa$ we have

$$\left\langle \left| c - \sum_{k=1}^{n} c_k \right|, \mathbb{Q}_\kappa \right\rangle \leq \sum_{k=n+1}^{\infty} \langle |c_k|, \mathbb{Q}_\kappa \rangle \leq \sum_{k=n+1}^{\infty} \|c_k\|_{\mathcal{M}}.$$

By taking the supremum over all $\mathbb{Q}_\kappa \in \mathcal{M}_\kappa$, it follows that $c_0 \in \mathcal{V}^{\mathcal{M}}$ and $\sum_{k=1}^{\infty} c_k = c_0$ in $\|\cdot\|_{\mathcal{M}}$. $\square$

REMARK 3.2. A norm of the form $\|\cdot\|_{\mathcal{M}}$ first appeared in Delbaen and Schachermayer (1997), who studied the Banach-space properties of the space of workable contingent claims.

At this point, we can introduce the third (and final) update of the notation of (3.2). Let $\mathcal{V}^{\mathcal{M}}_{\kappa+}$ denote the set of nonnegative elements in $\mathcal{V}^{\mathcal{M}}_\kappa$. For $c \in \mathcal{V}^{\mathcal{M}}_{\kappa+}$ a constant $y > 0$ and $\mathbb{Q} \in \mathcal{D}_\kappa(y)$, we define

$$(3.3) \qquad \langle c, \mathbb{Q} \rangle \triangleq \{\langle c', \mathbb{Q} \rangle : c' \in \mathbb{L}^{\infty}(\mathcal{O}, \mathbb{P}_\kappa)_+, \ c' \leq c \ \mathbb{P}_\kappa\text{-a.s.}\}.$$

Proposition 3.2 implies that $\langle c, \mathbb{Q} \rangle \leq y \|c\|_{\mathcal{M}} < \infty$ for any $\mathbb{Q} \in \mathcal{D}_\kappa(y)$. We can therefore extend the mapping $\langle \cdot, \cdot \rangle$ to a pairing (a bilinear form) between the vector spaces $\mathcal{V}^{\mathcal{M}}_\kappa$ and $\mathbf{ba}^{\mathcal{M}}$, where $\mathbf{ba}^{\mathcal{M}}$ is defined as the linear space spanned by $\mathcal{D}_\kappa$, that is,

$$\mathbf{ba}^{\mathcal{M}} \triangleq \{\mathbb{Q} \in \mathbf{ba}(\mathcal{O}, \mathbb{P}_\kappa) : \exists\, y > 0, \mathbb{Q}^+, \mathbb{Q}^- \in \mathcal{D}_\kappa(y) \text{ such that } \mathbb{Q} = \mathbb{Q}^+ - \mathbb{Q}^-\}.$$

The linear space $\mathbf{ba}^{\mathcal{M}}$ plays the role of the ambient space in which the dual domain is situated. It replaces the space $\mathbf{ba}$ appearing in Cvitanić, Schachermayer and Wang (2001) and Karatzas and Žitković (2003), and allows us to deal with unbounded random endowment and the stochastic clock.

In this way the action $\langle \cdot, \mathbb{Q} \rangle$ defined in (3.3) identifies $\mathbb{Q} \in \mathbf{ba}^{\mathcal{M}}$ with a linear functional on $(\mathcal{V}^{\mathcal{M}}, \|\cdot\|_{\mathcal{M}})$ and, by the construction of the pairing $\langle \cdot, \cdot \rangle$, the dual norm

$$\|\mathbb{Q}\|_{\mathbf{ba}^{\mathcal{M}}} \triangleq \sup_{c \in \mathcal{V}^{\mathcal{M}}_\kappa : \|c\|_{\mathcal{M}} \leq 1} |\langle c, \mathbb{Q} \rangle|$$

of $\mathbb{Q} \in \mathcal{D}_\kappa(y)$ (seen as a linear functional on $\mathcal{V}^{\mathcal{M}}_\kappa$) is at most equal to $2y$. We can, therefore, identify $\mathbf{ba}^{\mathcal{M}}$ with a subspace of the topological dual of $\mathcal{V}^{\mathcal{M}}_\kappa$ and $\mathcal{D}_\kappa(y)$ with its bounded subset. Moreover, by virtue of its definition as a polar set of $(\mathcal{M}_\kappa)^\circ$, $\mathcal{D}_\kappa(y)$ is closed in $\mathbf{ba}^{\mathcal{M}}$ in the $\sigma(\mathbf{ba}^{\mathcal{M}}, \mathcal{V}^{\mathcal{M}}_\kappa)$ topology, so that the following proposition becomes a direct consequence of Alaoglu's theorem.

PROPOSITION 3.4. *For every $y > 0$, $\mathcal{D}_\kappa(y)$ is $\sigma(\mathbf{ba}^{\mathcal{M}}, \mathcal{V}^{\mathcal{M}}_\kappa)$ compact.*

Finally, we state a version of the budget-constraint characterization of admissible consumption processes, rewritten to achieve a closer match with our newly introduced setup. It follows directly from Propositions 2.2 and 3.2.

PROPOSITION 3.5. *For any $y > 0$, $x \in \mathbb{R}$ and a nonnegative optional process $(c_t)_{t \in [0,T]}$, we have the equivalence*

$$c \in \mathcal{A}(x, \mathcal{E}) \iff y\langle c, \mathbb{Q}\rangle \leq xy + \langle e, \mathbb{Q}\rangle \quad \text{for all } \mathbb{Q} \in \mathcal{D}_\kappa(y),$$

*where $\mathcal{E}_t = \int_0^t e_u \, d\kappa_u$. Moreover, to check whether $c \in \mathcal{A}(x, \mathcal{E})$, it is enough to show $y\langle c, \mathbb{Q}\rangle \leq xy + \langle e, \mathbb{Q}\rangle$ for all $\mathbb{Q} \in \mathcal{M}_\kappa(y)$ only.*

## 4. The dual optimization problem and the main result.

4.1. *The convex conjugate $\mathbf{V}$ and related functionals.* We define a convex functional $\mathbf{V} : \mathbf{ba}^{\mathcal{M}} \to (-\infty, \infty]$ by

$$(4.1) \qquad \mathbf{V}(\mathbb{Q}) \triangleq \sup_{c \in \mathcal{V}_+^{\mathcal{M}}} (\mathbf{U}(c) - \langle c, \mathbb{Q}\rangle)$$

and call it the convex conjugate of $\mathbf{V}$. The functional $\mathbf{V}$ plays the central role in the convex-duality treatment of our utility-maximization problem.

By strict concavity and continuous differentiability of the mapping $x \mapsto U(\omega, t, x)$, there exists a unique random field $I : \Omega \times [0, T] \times (0, \infty)$ that solves the equation $U_x(\omega, t, I(\omega, t, y)) = y$. Using the random field $I$, we introduce a functional $\mathbf{I}$, defined on and taking values in the set of strictly positive optional process, by $\mathbf{I}(Y)_t(\omega) = I(\omega, t, Y_t)$. The functional $\mathbf{I}$ is called the inverse marginal utility functional. We note for the future use the well-known relationship

$$(4.2) \qquad \begin{aligned} U(\omega, t, I(\omega, t, y)) &= V(\omega, t, y) + yI(\omega, t, y), \\ &(\omega, t, y) \in \Omega \times [0, T] \times (0, \infty), \end{aligned}$$

where $V$ is the convex conjugate of the utility random field $U$, defined by $V(\omega, t, y) \triangleq \sup_{x>0}[U(\omega, t, x) - xy]$ for $(\omega, t, y) \in \Omega \times [0, T] \times (0, \infty)$.

For a function $f : X \to \overline{\mathbb{R}}$ with an arbitrary domain $X$, taking values in the extended set of real numbers $\overline{\mathbb{R}} = [-\infty, \infty]$, we adopt the standard notation $\mathrm{Dom}(f) = \{x \in X : f(x) \in (-\infty, \infty)\}$.

The following proposition represents the convex conjugate $\mathbf{V}$ in terms of the regular part of its argument, relating the definition (4.1) to the corresponding formulations in Cvitanić, Schachermayer and Wang (2001) and Karatzas and Žitković (2003).

PROPOSITION 4.1. *The domain* $\mathrm{Dom}(\mathbf{V})$ *of the convex conjugate* $\mathbf{V}$ *of* $\mathbf{U}$ *satisfies* $\mathrm{Dom}(\mathbf{V}) \subseteq \mathbf{ba}_+^{\mathcal{M}}$ *and* $\mathrm{Dom}(\mathbf{V}) + \mathbf{ba}_+^{\mathcal{M}} \subseteq \mathrm{Dom}(\mathbf{V})$. *For* $\mathbb{Q} \in \mathrm{Dom}(\mathbf{V})$, *we have* $\mathbf{V}(\mathbb{Q}) = \mathbf{V}(\mathbb{Q}^r)$, *where* $\mathbb{Q}^r \in \mathbf{ba}_+^{\mathcal{M}}$ *is the regular part of the finitely additive measure* $\mathbb{Q}$. *Moreover, there exists a nonnegative optional process* $Y^{\mathbb{Q}}$, *such that*

$$(4.3) \qquad \mathbf{V}(\mathbb{Q}) = \mathbb{E} \int_0^T V(t, Y_t^{\mathbb{Q}}) \, d\kappa_t.$$

*When* $\mathbb{Q}$ *is countably additive, the process* $(Y_t^{\mathbb{Q}})_{t \in [0,T]}$ *coincides with the synonymous martingale defined in Proposition* 3.1.

PROOF. For $\mathbb{Q} \notin \mathbf{ba}_+^{\mathcal{M}}$, there exists an optional set $A$ such that $q \triangleq -\mathbb{Q}[A] > 0$. For a constant $\varepsilon > 0$, we define a sequence $(c^n)_{n \in \mathbb{N}}$ of optional processes by $c^n \triangleq \varepsilon + n\mathbf{1}_A$. Let $G$ being the constant from part 1(b) of Definition 2.3. Then

$$\mathbf{V}(\mathbb{Q}) \geq \mathbf{U}(c^n) - \langle c^n, \mathbb{Q} \rangle \geq \mathbb{E} \int_0^T U(\omega, t, \varepsilon) \, d\kappa_t - \varepsilon + nq \geq G - \varepsilon + nq \to \infty$$

yields $\mathbf{V}(\mathbb{Q}) = \infty$ and so $\mathrm{Dom}(\mathbf{V}) \subseteq \mathbf{ba}_+^{\mathcal{M}}$. To show that $\mathrm{Dom}(\mathbf{V}) + \mathbf{ba}_+^{\mathcal{M}} \subseteq \mathrm{Dom}(\mathbf{V})$, we need only to note that it follows directly from the monotonicity of $\mathbf{V}$.

For the second claim, let $\mathbb{Q} \in \mathbf{ba}_+^{\mathcal{M}}$ and let $\mathrm{Sing}(\mathbb{Q})$ denote the family of all optional sets $A \subseteq [0,T] \times \Omega$ such that $\mathbb{Q}^s(A) = 0$, where $\mathbb{Q}^s$ denotes the singular part of the finitely additive measure $\mathbb{Q}$. For $A \in \mathrm{Sing}(\mathbb{Q})$, $\delta > 0$ and an arbitrary $c \in \mathcal{V}_+^{\mathcal{M}}$, we define an optional process $\hat{c} = \hat{c}^{(\delta, A)}$ by $\hat{c} \triangleq c\mathbf{1}_A + \delta c \mathbf{1}_{A^c}$. Excluding the trivial cases when $\mathbf{U}(c) = -\infty$ or $\mathbf{U}(c) = +\infty$, we assume $\mathbf{U}(c) \in \mathbb{R}$, so that Assumption 2.4 implies that $\mathbf{U}(\delta c), \mathbf{U}(\hat{c}) \in \mathbb{R}$, as well. Now

$$(4.4) \quad \begin{aligned} &\mathbf{U}(c) - \langle c, \mathbb{Q}^r \rangle - \mathbf{U}(\hat{c}) + \langle \hat{c}, \mathbb{Q} \rangle \\ &= \mathbb{E} \int_0^T (U(t, c_t) - U(t, \delta c_t)) \mathbf{1}_{A^c} \, d\kappa_t - (1-\delta) \langle c \mathbf{1}_{A^c}, \mathbb{Q}^r \rangle + \delta \langle c, \mathbb{Q}^s \rangle. \end{aligned}$$

According to Bhaskara Rao and Bhaskara Rao [(1983), Theorem 10.3.2, page 234], $\mathrm{Sing}(\mathbb{Q}_\kappa)$ contains sets with the $\mathbb{P}_\kappa$ probability arbitrarily close to 1, so we can make the right-hand side of the expression in (4.4) arbitrarily small in absolute value by a suitable choice of $A \in \mathrm{Sing}(\mathbb{Q})$ and $\delta$. It follows immediately that

$$\mathbf{V}(\mathbb{Q}^r) = \sup_{c \in \mathcal{V}^{\mathcal{M}}} [\mathbf{U}(c) - \langle c, \mathbb{Q}^r \rangle] \leq \sup_{c \in \mathcal{V}^{\mathcal{M}}} [\mathbf{U}(c) - \langle c, \mathbb{Q} \rangle] = \mathbf{V}(\mathbb{Q})$$

and the equality $\mathbf{V}(\mathbb{Q}) = \mathbf{V}(\mathbb{Q}^r)$ follows from the monotonicity of $\mathbf{V}$.

Note further that $\mathbb{Q}^r$ is a countably additive measure on the $\sigma$-algebra of optional sets, absolutely continuous with respect to the measure $\mathbb{P}_\kappa$. It follows by the Radon–Nikodym theorem that the optional process $(Y_t^{\mathbb{Q}})_{t \in [0,T]}$ defined by

$$(4.5) \qquad Y^{\mathbb{Q}}(t, \omega) = \frac{d\mathbb{Q}^r}{d\mathbb{P}_\kappa} \quad \text{satisfies} \quad \langle c, \mathbb{Q}^r \rangle = \mathbb{E} \int_0^T c_t Y_t^{\mathbb{Q}} \, d\kappa_t.$$

Let us now combine the representation (4.5) with the fact that $\mathbf{V}(\mathbb{Q}) = \mathbf{V}(\mathbb{Q}^r)$. By the definition of the convex conjugate function $V$,

$$\mathbf{V}(\mathbb{Q}) = \mathbf{V}(\mathbb{Q}^r) = \sup_{c \in \mathcal{V}_+^{\mathcal{M}}} (\mathbf{U}(c) - \langle c, \mathbb{Q}^r \rangle)$$

$$= \sup_{c \in \mathcal{V}_+^{\mathcal{M}}} \mathbb{E} \int_0^T (U(t, c(t)) - c(t) Y_t^{\mathbb{Q}}) \, d\kappa_t \leq \mathbb{E} \int_0^T V(t, Y_t^{\mathbb{Q}}) \, d\kappa_t.$$

The reverse inequality follows from the differentiability of the function $V(t, \cdot)$ by taking a bounded sequence in $\mathcal{V}^{\mathcal{M}}$, which converges to $-\frac{\partial}{\partial y} V(t, y)$ monotonically, in the supremum that defines $\mathbf{V}(\mathbb{Q}^r)$. $\square$

REMARK 4.1. The action of the functional $\mathbf{I}$ can be extended to the set of all $\mathbb{Q} \in \mathbf{ba}_+^{\mathcal{M}}$ that satisfy $Y_t^{\mathbb{Q}} > 0$ $\mathbb{P}_\kappa$-a.e. by $\mathbf{I}(\mathbb{Q})_t \triangleq \mathbf{I}(Y^{\mathbb{Q}})_t$, obtaining immediately $\mathbf{I}(\mathbb{Q}) = \mathbf{I}(\mathbb{Q}^r)$.

4.2. *The dual problem.* The convex conjugate $\mathbf{V}$ serves as the main ingredient in the convex-duality treatment of the primal problem. We start by introducing the dual problem, with the value function $v$:

$$v(y) \triangleq \inf_{\mathbb{Q} \in \mathcal{D}_\kappa(y)} \mathbf{V}^{\mathcal{E}}(\mathbb{Q}),$$

(4.6)
$$y \in [0, \infty) \text{ where } \mathbf{V}^{\mathcal{E}}(\mathbb{Q}) \triangleq \mathbf{V}(\mathbb{Q}) + \langle e, \mathbb{Q} \rangle.$$

For $y < 0$, we set $v(y) = +\infty$ and note that $v(0) < \infty$ precisely when the utility functional $\mathbf{U}$ is bounded from above.

4.3. *The main result.* Finally we state our central result in the following theorem. The proof is given through a number of auxiliary results in the Appendix A.

THEOREM 4.2. *Let the financial market* $(S_t^i)_{t \in [0,T]}$, $i = 1, \ldots, d$, *be arbitrage-free as in Assumption* 2.1 *and let the random endowment process* $(\mathcal{E}_t)_{t \in [0,T]}$ *admit a density* $(e_t)_{t \in [0,T]}$ *so that* $\mathcal{E}_t = \int_0^t e_u \, d\kappa_u$, *where* $(\kappa_t)_{t \in [0,T]} \in \mathcal{V}_+$ *is a stochastic clock. Let* $U$ *be a utility random field as defined in Definition* 2.3 *and let* $\mathbf{U}$ *be the corresponding utility functional. If* $\mathbf{U}$ *satisfies Assumption* 2.4 *and the value function* $u$ *satisfies Assumption* 2.9, *then:*

1. *The concave value function $u(\cdot)$ is finite and strictly increasing on $(-\mathcal{L}(\mathcal{E}), \infty)$ and $u(x) = -\infty$ for $x < -\mathcal{L}(\mathcal{E})$, where $\mathcal{L}(\mathcal{E}) \triangleq \inf_{\mathbb{Q} \in \mathcal{M}} \mathbb{E}^{\mathbb{Q}}[\mathcal{E}_T]$ denotes the lower hedging price of the contingent claim $\mathcal{E}_T$.*
2. *We have $\lim_{x \to (-\mathcal{L}(\mathcal{E}))+} u'(x) = +\infty$ and $\lim_{x \to \infty} u'(x) = 0$.*
3. *The dual value function $v(\cdot)$ is finitely valued and continuously differentiable on $(0, \infty)$ and $v(y) = +\infty$ for $y < 0$.*
4. *We have $\lim_{y \to 0+} v'(y) = -\infty$ and $\lim_{y \to \infty} v'(y) = -\mathcal{L}(\mathcal{E})$.*
5. *For any $y \geq 0$, there exists a solution to the dual problem (4.6), that is, $v(y) = \mathbf{V}(\widehat{\mathbb{Q}}^y) + \langle e, \widehat{\mathbb{Q}}^y \rangle$ for some $\widehat{\mathbb{Q}}^y \in \mathcal{D}_\kappa(y)$.*
6. *For $x > -\mathcal{L}(\mathcal{E})$, the primal problem has a solution $(\hat{c}_t^x)_{t \in [0,T]}$ that is unique $d\kappa$-a.e.*
7. *The unique solution $(\hat{c}_t^x)_{t \in [0,T]}$ of the primal problem is of the form $\hat{c}_t^x = \mathbf{I}(\widehat{\mathbb{Q}}^y)_t$, where $\widehat{\mathbb{Q}}^y$ is a solution of the dual problem that corresponds to $y > 0$ such that $x = -v'(y)$.*

4.4. *A closer look at the dual domain.* Given that the solution of the primal problem can be expressed as a function of the process $(Y_t^{\mathbb{Q}})_{t \in [0,T]}$ from Proposition 4.1, it is useful to have more information on its probabilistic structure. When $\mathbb{Q} \in \mathcal{M}_\kappa$, Proposition 3.1 implies that $Y^{\mathbb{Q}}$ is a nonnegative cadlag martingale. In general, we can only establish the supermartingale property for a (large enough) subclass of ($\mathbb{P}_\kappa$-a.s.) maximal processes in $\{Y^{\mathbb{Q}} : \mathbb{Q} \in \mathcal{D}(1)\}$. In the contrast with the case studied in Karatzas and Žitković (2003), we cannot establish any strong trajectory regularity properties such as right continuity and have to satisfy ourselves with the weaker property of optional measurability.

PROPOSITION 4.3. *For $\mathbb{Q} \in \mathcal{D}(1)$ there exists an optional process $(F_t)_{t \in [0,T]}$, taking values in $[0, 1]$, and $\mathbb{Q}' \in \mathcal{D}(1)$ such that the following statements hold:*

1. *We have $Y_t^{\mathbb{Q}} = Y_t^{\mathbb{Q}'} F_t$.*
2. *The process $(Y_t^{\mathbb{Q}'})_{t \in [0,T]}$ has a $d\kappa$ version which is an optional supermartingale.*
3. *There exists a sequence of martingale measures $\{\mathbb{Q}_n\}_{n \in \mathbb{N}}$ such that $Y^{\mathbb{Q}_n} \to Y^{\mathbb{Q}'}$ $d\kappa$-a.e.*

PROOF. We start by observing that $\mathbb{E}[\int_0^T Y_t^{\mathbb{Q}} c(t) \, d\kappa_t] \leq \langle c, \mathbb{Q} \rangle \leq 1$ for all $c \in \mathcal{A}(1, 0)$. In other words, $Y^{\mathbb{Q}}$ is in the $\mathbb{P}_\kappa$ polar set of $\mathcal{A}(1, 0)$ in the terminology of Brannath and Schachermayer (1999). By the characterization in Proposition 3.5, $\mathcal{A}(1, 0)$ can be written as the polar of $\mathcal{M}_\kappa$, and the bipolar theorem of Brannath and Schachermayer (1999) states that $Y^{\mathbb{Q}}$ is an element of the smallest convex, solid and closed (in $\mathbb{P}_\kappa$ probability) set containing $\mathcal{M}_\kappa$. Therefore, there exists a process $(F_t)_{t \in [0,T]}$, taking values in $[0,1]$, and an optional process $(Y_t)_{t\in[0,T]}$, ($\mathbb{P}_\kappa$-a.s.) maximal in the bipolar of $\mathcal{M}_\kappa$, such that $Y_t^{\mathbb{Q}} = Y_t F_t$. Moreover, the same theorem implies that there exists a sequence $\{\mathbb{Q}^{(n)}\}_{n\in\mathbb{N}}$ in $\mathcal{M}$ and a sequence $\{F^{(n)}\}_{n\in\mathbb{N}}$ of optional processes taking values in $[0,1]$, such that $Y_t^{\mathbb{Q}^{(n)}} F_t^{(n)} \to Y_t$ $\mathbb{P}_\kappa$ a.s. The sequence of positive processes $Y^{\mathbb{Q}^{(n)}}$ is bounded in $\mathbb{L}^1(\mathbb{P}_\kappa)$; thus the theorem of Komlós [see Schwartz (1986)] asserts the existence of a nonnegative optional process $(\widetilde{Y}_t)_{t\in[0,T]}$ and a sequence of finite convex combinations of the elements of the sequence $\{\mathbb{Q}^{(n)}\}_{n\in\mathbb{N}}$ (still denoted by $\{\mathbb{Q}^{(n)}\}_{n\in\mathbb{N}}$) such that $Y_t^{\mathbb{Q}^{(n)}} \to \widetilde{Y}_t$ $\mathbb{P}_\kappa$-a.s. It is now a simple consequence of Fatou's lemma that $\widetilde{Y}$ is an element of the bipolar of $\mathcal{M}_\kappa$ dominating $Y_t$. Since $Y_t$ is maximal, we conclude that $\widetilde{Y}_t = Y_t$ $\mathbb{P}_\kappa$-a.s. The supermartingale property of $(Y)_{t\in[0,T]}$ follows from Fatou's lemma applied to the sequence $\{(Y_t^{\mathbb{Q}^{(n)}})_{t\in[0,T]}\}_{n\in\mathbb{N}}$.

We are left now with the task of producing $\mathbb{Q}' \in \mathcal{D}(1)$ such that $Y_t = Y_t^{\mathbb{Q}'}$. To do that, take $\mathbb{Q}'$ to be any cluster point of the sequence $\{\mathbb{Q}^{(n)}\}_{n\in\mathbb{N}}$ in $\mathcal{D}(1)$ in the $\sigma(\mathbf{ba}^{\mathcal{M}}, \mathcal{V}_\kappa^{\mathcal{M}})$ topology. The existence of such a $\mathbb{Q}'$ is guaranteed by Proposition 3.4. Finally, it is a consequence of Cvitanić, Schachermayer and Wang [(2001), Lemma A.1, page 16] that $Y_t = Y_t^{\mathbb{Q}'}$ $\mathbb{P}_\kappa$-a.s. □

**5. An example.** To illustrate the theory developed so far, in this section we present an example of a utility-maximization problem with a random clock given by the local time at 0 of an Ornstein–Uhlenbeck process.

5.1. *Description of the market model.* Let $(B_t, W_t)_{t\in[0,\infty)}$ be two correlated Brownian motions defined on a probability space $(\Omega, \mathcal{F}, \mathbb{P})$ and let $(\mathcal{F}_t)_{t\in[0,\infty)}$ be the filtration they generate, augmented by the $\mathbb{P}$-null sets to satisfy the usual conditions. We assume that the correlation coefficient $\rho \in (-1,1)$ is fixed so that $d[B,W]_t = \rho \, dt$.

The financial market consists of one riskless asset $S_t^0 \equiv 1$ and a risky asset $(S_t)_{t\in[0,\infty)}$ which satisfies

$$dS_t = S_t(\mu \, dt + \sigma \, dB_t), \qquad S_0 = s_0,$$

where $\mu \in \mathbb{R}$ is the stock appreciation rate and $\sigma > 0$ is the volatility.

Apart from the tradeable asset $(S_t)_{t\in[0,\infty)}$, there is an Orstein–Uhlenbeck process $(R_t)_{t\in[0,\infty)}$ defined as the unique strong solution of

$$dR_t = -\alpha R_t \, dt + dW_t, \qquad R_0 = 0.$$

We call $(R_t)_{t\in[0,\infty)}$ the *index process* and interpret it as the process that models a certain state variable of the economy, possibly related to political stability or some aspect of the government's economic policy. The index process is nontradable and its role is to impose constraints on the consumption: We are allowed to withdraw money from the trading account only when $|R_t| < \varepsilon$. An agent with an initial endowment $x$ and a utility random field $U(\cdot, \cdot, \cdot)$ then naturally tries to choose a strategy so as to maximize the utility of consumption of the form

$$\mathbb{E} \int_0^\tau U(\omega, t, c(t)) \mathbf{1}_{\{|R_t| < \varepsilon\}} \, dt \tag{5.1}$$

on some trading horizon $[0, \tau]$. If we introduce the notation $\kappa_t^\varepsilon = \frac{1}{\varepsilon} \int_0^t \mathbf{1}_{\{|R_t| < \varepsilon\}} \, dt$, the expression in (5.1) becomes (up to a multiplicative constant)

$$\mathbb{E} \int_0^\tau U(\omega, t, c(t)) \, d\kappa_t^\varepsilon. \tag{5.2}$$

Assuming that $\varepsilon$ is a small constant, the process $\kappa^\varepsilon$ can be approximated by the local time $\kappa_t$ of the process $R_t$. We define the time horizon $\tau = \tau_1$, where $\tau_s \triangleq \inf\{t > 0 : \kappa_t > s\}$ is the inverse local-time process. In this way our agent gets exactly one unit of consumption time (as measured by the clock $\kappa$) from the start to the end of the trading interval. It is, therefore, our goal to solve the following problem, defined in terms of its value function $u(\cdot)$:

$$u(x) = \sup_{c \in \mathcal{A}(x,0)} \mathbb{E} \int_0^{\tau_1} U(\omega, t, c_t) \, d\kappa_t, \qquad x > 0. \tag{5.3}$$

5.2. *Absence of arbitrage.* The time horizon $\tau$ defined above is clearly not a bounded random variable, so the results in the main body of this paper do not apply directly. However, to pass from an infinite to a finite horizon, it is enough to apply a deterministic time change that maps $[0, \infty)$ onto $[0, 1)$ and to note that no important part of the structure of the problem is lost in this way (we leave the easy details of the argument to the reader). Of course, we need to show that all the assumptions of Theorem 4.2 are satisfied. The validity of Assumption 2.9 has to be checked on a case-by-case basis (see Remark 5.1 for the case of log utility). Therefore, we are left with Assumption 2.1. To proceed, we need to exhibit a countably additive probability measure $\mathbb{Q}$ equivalent to $\mathbb{P}$ such that the asset-price process $(S_t)_{t \in [0, \infty)}$ is a $\mathbb{Q}$-local martingale on the stochastic interval $[0, \tau_1]$. The obvious candidate is the measure $\mathbb{Q}_0$ defined in terms of its Radon–Nikodym derivative with respect to $\mathbb{P}$ by

$$\frac{d\mathbb{Q}_0}{d\mathbb{P}} = Z_{\tau_1}^0 \qquad \text{where } Z_{\tau_1}^0 \triangleq \exp(-\theta B_{\tau_1} - \tfrac{1}{2}\theta^2 \tau_1) \tag{5.4}$$

and $\theta = \mu/\sigma$ is the market price of risk coefficient. Once we show that $\mathbb{E}[Z_{\tau_1}^0] = 1$, it follows directly from Girsanov's theorem [see Karatzas and Shreve (1991), Theorem 3.5.1, page 191] that $(S)_{t \in [0, \infty)}$ is a $\mathbb{Q}$-local martingale on $[0, \tau_1]$. The equivalence of the measures $\mathbb{Q}_0$ and $\mathbb{P}$ is a consequence of the fact that $\tau_1 < \infty$ a.s, which follows from the following proposition which lists some distributional properties of the process $(R_t)_{t\in[0,\infty)}$ and its local time $(\kappa_t)_{t\in[0,\infty)}$.

PROPOSITION 5.1. *For $\xi < 0$ and $x \geq 0$, let $H_\xi(x)$ denote the value of the Hermite function*

$$(5.5) \qquad H_\xi(x) = \frac{1}{2\Gamma(-\xi)} \int_0^\infty e^{-s-2x\sqrt{s}} s^{-\xi/2-1}\,ds.$$

*For the Ornstein–Uhlenbeck process $(R_t)_{t\in[0,\infty)}$ and the inverse $(\tau_s)_{s\in[0,\infty)}$ of its local time at 0 $(\kappa_t)_{t\in[0,\infty)}$, we have the explicit expressions*

$$(5.6) \qquad \mathbb{E}[\exp(-\lambda\tau_s)|R_0=0] = \begin{cases} \exp(-s\psi(\lambda)), & \lambda > -\alpha, \\ \infty, & \lambda \leq -\alpha, \end{cases}$$

*where the Laplace exponent $\psi(\lambda)$ is given by*

$$(5.7) \qquad \psi(\lambda) = \alpha \frac{2^{1+\lambda/\alpha}\Gamma(1/2+\lambda/2\alpha)^2}{\sqrt{2\pi}\Gamma(\lambda/\alpha)},$$

*and, with $T_0 = \inf\{t>0 : R_t = 0\}$, we have*

$$(5.8) \qquad \mathbb{E}[\exp(-\lambda T_0)|R_0=r] = j(\lambda,|r|),$$

*where*

$$j(\lambda,r) \triangleq 2^{\lambda/\alpha}\frac{\Gamma((1+\lambda/\alpha)/2)}{\Gamma(1/2)} H_{-\lambda/\alpha}\left(\frac{r}{\sqrt{2}}\right).$$

PROOF. See equation (2.0.1) of Borodin and Salminen [(2002), page 542] for (5.6) and equation (4.0.1) of Borodin and Salminen [(2002), page 557] for (5.8). Use the identity $D_\zeta(x) = 2^{-\zeta/2}\exp(-x^2/4)H_\zeta(x/\sqrt{2})$. □

To prove the equality $\mathbb{E}[Z^0_{\tau_1}] = 1$, it is be enough to show that $\mathbb{E}[\exp(\frac{1}{2}\theta^2\tau_1)] < \infty$ by the Novikov's criterion [Karatzas and Shreve (1991), Proposition 3.5.12, page 198]. Equation (5.6) of Proposition 5.1 implies that for $\alpha > \theta^2/2$, we have $\mathbb{E}[\exp(\frac{1}{2}\theta^2\tau_1)] < \infty$, which proves the following proposition.

PROPOSITION 5.2. *When $\alpha > \theta^2/2$, there is no arbitrage on the stochastic interval $[0,\tau_1]$.*

5.3. *The optimal consumption and portfolio choice.* It was shown in Karatzas and Žitković (2003) that the maximal dual processes in the context of the financial markets driven by Itô processes with bounded coefficients are in fact local martingales and their structure was described. This result can be extended to our case as follows.

THEOREM 5.3. *Let the utility random field $U$ satisfy Assumptions* 2.4 *and* 2.9. *Then, for $x > 0$, there exists a predictable process $(\nu_t^x)_{t \in [0,\infty)}$ such that the $\mathbb{P}_\kappa$-a.e. unique solution $(\hat{c}_t^x)_{t \in [0,\infty)}$ of the problem posed in* (5.3) *is given by $\hat{c}_t^x(\omega) = I(\omega, t, Z_t^{\nu^x}(\omega))$. The process $(Z_t^{\nu^x})_{t \in [0,\infty)}$ is a local martingale that satisfies*

$$(5.9) \qquad dZ_t^{\nu^x} = Z_t^{\nu^x}(\nu_t^x \, dW_t - (\theta + \rho \nu_t^x) \, dB_t), \qquad Z_0^{\nu^x} = y,$$

*where $y > 0$ is the unique solution of $-v'(y) = x$. The portfolio process $(\pi_t^x)_{t \in [0,\infty)}$ that finances $(\hat{c}_t^x)_{t \in [0,\infty)}$ and the process $(\nu_t^x)_{t \in [0,\infty)}$ are given by*

$$(5.10) \qquad \pi_t^x = \frac{X_t}{\sigma S_t}(\theta + \rho \nu_t^x) + \frac{\psi_t^B}{\sigma S_t Z_t^{\nu^x}}, \qquad \nu_t^x = \frac{1}{X_t Z_t^{\nu^x}} \psi_t^W,$$

*where $(X_t)_{t \in [0,\infty)}$ is the wealth process that corresponds to $(\pi_t^x)_{t \in [0,\infty)}$ and $(\hat{c}_t^x)_{t \in [0,\infty)}$, given by*

$$(5.11) \qquad dX_t = \pi_t^x \, dS_t - \hat{c}_t^x \, d\kappa_t, \qquad X_0 = x,$$

*and $(\psi^B)_{t \in [0,\infty)}$ and $(\psi^W)_{t \in [0,\infty)}$ are predictable processes such that*

$$(5.12) \qquad xy + \int_0^{\tau_1} \psi_t^B \, dB_t + \int_0^{\tau_1} \psi_t^W \, dW_t = \int_0^{\tau_1} Z_t^{\nu^x} \hat{c}_t^x \, d\kappa_t.$$

PROOF. By Theorem 4.2, there exists a $\mathbb{P}_\kappa$-a.e. unique optimal consumption density $\hat{c}^x \in \mathcal{A}(x, 0)$ given by $\hat{c}_t^x = I(t, Y_t^{\mathbb{Q}})$ for some $\mathbb{Q} \in \mathcal{D}_\kappa(y)$. Since $(Y_t^{\mathbb{Q}})_{t \in [0,\infty)}$ solves the dual optimization problem and is, therefore, $\mathbb{P}_\kappa$-a.e. maximal, Proposition 4.3 states that there exist a sequence $\{\mathbb{Q}^{(n)}\}_{n \in \mathbb{N}}$ in $\mathcal{M}$ such that $Y^{\mathbb{Q}^{(n)}} \to Y^{\mathbb{Q}}$ $\mathbb{P}_\kappa$-a.s. By taking a further sequence of convex combinations which exist thanks to Komlós's theorem [see Komlós (1967) and Schwartz (1986)], we can assume that $Y_T^{\mathbb{Q}^{(n)}} \to Y_T^{\mathbb{Q}}$ $\mathbb{P}$-a.s. and $Y_t^{\mathbb{Q}^{(n)}} \to Y_t^{\mathbb{Q}^{(n)}}$ $\mathbb{P} \times \lambda$-a.e. Without going into tedious but straightforward details, we note that it is the consequence of the continuity of local martingales on Brownian filtrations, the filtered bipolar theorem [Žitković (2002), Theorem 2], and Lemma 2.5, Theorem 2.10 and Proposition 4.1 in Karatzas and Žitković (2003) that $(Y_t^{\mathbb{Q}})_{t \in [0,\infty)}$ possesses a $\mathbb{P}_\kappa$ version of the form $Y_t^{\mathbb{Q}} = y Z_t^\nu$, where $Z^\nu$ is a local martingale of the form (5.9).

Knowing that $\hat{c}^x \in \mathcal{A}(x, 0)$, there exists a portfolio process $(\pi_t^x)_{t \in [0,\infty)}$ such that the wealth process $(X_t)_{t \in [0,\infty)}$ given by (5.11) satisfies $X_{\tau_1} \geq 0$. The saturation of the budget constraint (see Lemma A.3.2) forces $X_{\tau_1} = 0$. Itô's lemma shows that the process

$$(5.13) \qquad M_t = X_t Z_t^\nu + \int_0^t Z_u^\nu \hat{c}_u^x \, d\kappa_u$$

is a nonnegative local martingale with $M_{\tau_1} = \int_0^{\tau_1} Z_u^\nu \hat{c}_u^x \, d\kappa_u$. By Lemma A.3.2, we have $\mathbb{E}[M_{\tau_1}] = x = M_0$. Therefore, $M$ is a martingale on $[0, \tau_1]$. The second equality in (5.10) follows by applying Itô's formula to (5.13) and equating coefficients with those in the expansion (5.12). $\square$

5.4. *The case of logarithmic utility.* To get explicit results, we consider now the agent whose utility function has the form $U(\omega, t, x) = \exp(-\beta t)\log(x)$, where the impatience rate $\beta$ is a positive constant. The expressions (5.10) prove indispensable because it is possible to get an explicit expression for the processes $(\psi_t^W)_{t \in [0,\infty)}$ and $(\psi_t^B)_{t \in [0,\infty)}$ from (5.12). The key feature of the logarithmic utility that allows us to do this is the fact that the inverse marginal utility function $I$ is given by $I(t, y) = \exp(-\beta t)/y$, so that the right-hand side of (5.12) becomes

$$(5.14) \qquad M_{\tau_1} \triangleq \int_0^{\tau_1} Z_t^\nu \hat{c}_t^x \, d\kappa_t = \int_0^{\tau_1} e^{-\beta t} \, d\kappa_t.$$

To progress with the explicit representation of the processes $(\psi_t^W)_{t \in [0,\infty)}$ and $(\psi_t^B)_{t \in [0,\infty)}$ from (5.12), in the following lemma we prove a useful fact about the conditional $\beta$ potential of the local time $(\kappa_t)_{t \in [0,\infty)}$, that is, the random process $(G_t)_{t \in [0,\infty)}$ defined by $G_t \triangleq \mathbb{E}[\int_0^{\tau_1} \exp(-\beta u) \, d\kappa_u | \mathcal{F}_t]$.

LEMMA 5.4. *A version of the process $G$ is given by*

$$(5.15) \quad G_t = \begin{cases} \exp(-\beta t) j(\beta, |R_t|) \dfrac{1 - \exp(-(1-\kappa_t)\Psi(\beta))}{\Psi(\beta)} + \displaystyle\int_0^t e^{-\beta u} \, d\kappa_u, & \kappa_t \leq 1, \\ \displaystyle\int_0^{\tau_1} e^{-\beta u} \, d\kappa_u, & \kappa_t > 1, \end{cases}$$

*where the functions $\psi$ and $j$ are defined in (5.7) and (5.8).*

PROOF. We start by defining a family of stopping times $T_0(t) = \inf\{u \geq t : R_u = 0\}$ and note that because $d\kappa_u$ does not charge the complement of the zero set of $R_t$, we have

$$(5.16) \qquad G_t = \mathbb{E}\left[\int_{T_0(t)}^{\tau_1} e^{-\beta u} \, d\kappa_u \Big| \sigma(\kappa_t, R_t)\right] + \int_0^t e^{-\beta u} \, d\kappa_u.$$

Replacement of the $\sigma$-algebra $\mathcal{F}_t$ by $\sigma(\kappa_t, R_t)$ is permitted by the Markov property of the process $(\kappa_t, R_t)$.

When $\kappa_t \geq 1$, the value of $G_t$ is trivially given by (5.15), so we can restrict our attention to the value of the function $g(t, r, k) = \mathbb{E}[\int_{T_0(t)}^{\tau_1} e^{-\beta u} \, d\kappa_u | \kappa_t = k, R_t = r]$ for $k < 1$, because then (5.16) implies that $G_t = g(t, R_t, \kappa_t) + \int_0^t \exp(-\beta u) \, d\kappa_u$ on $\{\kappa_t < 1\}$. Using again the strong Markov property and time homogeneity of $(\kappa_t, R_t)$, we obtain

$$\begin{aligned}(5.17) \quad g(t, r, k) &= \mathbb{E}\left[e^{-\beta T_0(t)} \int_{T_0}^{\tau_1} e^{-\beta(u - T_0(t))} \, d\kappa_u \Big| R_t = r, \kappa_t = k\right] \\ &= e^{-\beta t} \mathbb{E}[e^{-\beta T_0(0)} | R_0 = r] \mathbb{E}\left[\int_0^{\tau_{1-k}} e^{-\beta t} d\kappa_t \Big| R_0 = 0, \kappa_0 = 0\right].\end{aligned}$$

The second term in the above expression is given in (5.8). As for the third term, a change of variables yields

$$(5.18) \quad \mathbb{E}\left[\int_0^{\tau_{1-k}} e^{-\beta t}\, d\kappa_t\right] = \int_0^{1-k} \mathbb{E}[e^{-\beta \tau_u}]\, du = \frac{1 - e^{-(1-k)\psi(\beta)}}{\psi(\beta)}. \qquad \Box$$

We have developed all the tools required to prove the following result

PROPOSITION 5.5. *In the setup of Theorem 5.3, set $U(\omega, t, x) = \exp(-\beta t) \times \log(x)$. Then we have the following explicit representations of the processes $(\pi_t^x)_{t \in [0, \infty)}$, $(\nu_t^x)_{t \in [0, \infty)}$ and $(\hat{c}_t^x)_{t \in [0, \infty)}$:*

$$(5.19) \quad \nu_t^x = -\mathrm{sgn}(R_t) h\left(\frac{|R_t|}{\sqrt{2}}\right) \qquad \text{where } h(z) \triangleq -\frac{2\beta}{\alpha} \frac{H_{-\beta/\alpha - 1}(z)}{H_{-\beta/\alpha}(z)},$$

$$(5.20) \quad \pi_t^x = \frac{X_t}{\sigma S_t}\left(\theta + \rho\, \mathrm{sgn}(R_t) h\left(\frac{|R_t|}{\sqrt{2}}\right)\right),$$

$$(5.21) \quad \hat{c}_t^x = X_t \frac{1 - \exp(-\Psi(\beta))}{(1 - \exp(-(1 - \kappa_t)\Psi(\beta)))}.$$

*Finally, the process $(\nu_t^x)_{t \in [0, T]}$ is bounded and so the optimal dual process $(Z_t^{\nu^x})_{t \in [0, T]}$ is a martingale.*

PROOF. Use of the Itô–Tanaka formula and expression (5.15) yields

$$(5.22) \quad \begin{aligned} \psi_t^B &= 0 \quad \text{and} \\ \psi_t^W &= \exp(-\beta t)\, \mathrm{sgn}(R_t) \frac{\partial}{\partial r} j(\beta, |R_t|) \frac{1 - \exp(-(1 - \kappa_t)\Psi(\beta))}{\Psi(\beta)}. \end{aligned}$$

Moreover, the martingale property of process $M_t$ from (5.13) implies that $X_t Z_t^{\nu^x} = G_t - \int_0^t e^{-\beta u}\, d\kappa_u$, and so (5.8), (5.10) and (5.12) can be combined into the explicit expression of the optimal dual process

$$\nu_t^y = \mathrm{sgn}(R_t) \frac{(\partial/\partial \beta) j(\beta, |R_t|)}{j(\beta, |R_t|)}.$$

Representation (5.8) and the identity $\frac{\partial}{\partial x} H_\xi(x) = 2\xi H_{\xi - 1}(x)$ [see Lebedev (1972), equation 10.5.2, page 289] complete the proof of (5.19).

Theorem 4.2 part 7 and identities (5.10) and (5.22) imply that

$$\hat{c}_t^x = \frac{X_t \Psi(\beta)}{y j(\beta, |R_t|)(1 - \exp(-(1 - \kappa_t)\Psi(\beta)))},$$

where $y$ satisfies $x = -v'(y)$. To get a more explicit expression for $y$, we combine (5.14) and (5.12) to get $xy = \mathbb{E}[\int_0^{\tau_1} \exp(-\beta t)\, d\kappa_t]$. After repeating the calculation in (5.18) with $k = 0$, we need only to rearrange the terms and remember that $R_t = 0$ $d\kappa$-a.e. to obtain (5.21).

We are left with the proof of the boundedness of the process $(\nu_t^x)_{t \in [0,\infty)}$. The asymptotic formula 10.6.3 in Lebedev [(1972), page 291] implies that $H_\xi(x) \sim C_\xi x^\xi$ as $x \to \infty$ for some positive constant $C_\xi$ depending on $\xi < 0$. Therefore, there exists a constant $D > 0$ such that $h(x) \sim Dx^{-1}$ as $x \to \infty$. Because of the existence of the limit $\lim_{x \to 0+} h(x)$, we conclude that $h$ is a bounded function on $[0, \infty)$. Hence, $(\nu_t^x)_{t \in [0,\infty)}$ is a bounded process, making $(Z_t^{\nu^x})_{t \in [0,T]}$ a martingale. $\square$

REMARK 5.1. In the generic setup of Theorem 5.3, we have explicitly assumed that $u(x) < \infty$ for at least one $x > 0$. In the case of the logarithmic utility random field treated above, the validity of such an assumption is implied by the chain of inequalities in which $\mathbb{Q}_0$ and $Z_{\tau_1}^0$ are as in (5.4), that is

$$u(x) - x = \sup_{c \in \mathcal{A}(x,0)} (\mathbf{U}(c) - x) \leq \mathbf{V}(\mathbb{Q}_0) = \mathbb{E}\int_0^{\tau_1}(-1 - \log(Z_t^0))\,d\kappa_t$$

(5.23)
$$\leq \mathbb{E}\bigg[\int_0^{\tau_1} \frac{1}{2}(\theta B_t^2 + 1 + \theta^2 t)\,d\kappa_t\bigg] = \frac{1}{2}\int_0^1 \mathbb{E}[\theta(1 + B_{\tau_s}^2) + \theta^2 \tau_s]\,ds$$

$$\leq \frac{\theta}{2} + \frac{(\theta^2 + 1)}{2}\int_0^1 \mathbb{E}[\tau_s]\,ds \leq \frac{\theta + (\theta^2 + 1)\mathbb{E}[\tau_1]}{2} < \infty.$$

The fact that $\mathbb{E}[\tau_1] < \infty$ [which can easily be deduced from (5.6)] implies both the final inequality in (5.23) and the equality $\mathbb{E}[B_{\tau_1}^2] = \mathbb{E}[\tau_1]$ through Wald's identity [see Problem 2.12, page 141 in Karatzas and Shreve (1991)].

## APPENDIX: A CONVEX-DUALITY PROOF OF THEOREM 4.2

We have divided the proof into several steps, each of which is stated as a separate lemma. Throughout this section all the conditions of Theorem 4.2 are assumed to be satisfied.

LEMMA A.1 (Global properties of the value functions). *The value function $u(\cdot)$ is convex, nondecreasing and $[-\infty, \infty)$ valued, while $v$ is concave and $(-\infty, \infty]$ valued. Moreover, the primal and the dual value functions $u(\cdot)$ and $v(\cdot)$ are convex conjugates of each other.*

PROOF. 1. Concavity of $u(\cdot)$ and convexity of $v(\cdot)$ are inherited from the properties of the objective functions $\mathbf{U}(\cdot)$ and $\mathbf{V}(\cdot)$ [see Ekeland and Témam (1999), proof of Lemma 2.1, page 50, for the standard argument]. The increase of $u(\cdot)$ follows from the inclusion $\mathcal{A}(x, \mathcal{E}) \subseteq \mathcal{A}(x', \mathcal{E})$ for $x < x'$.

2. By the Assumption 2.9, there exists $\tilde{x} \in \mathbb{R}$ such that $u(\tilde{x}) < \infty$. It follows immediately by concavity of $u(\cdot)$ that $u(x) < \infty$ for all $x \in \mathbb{R}$.

3. To establish the claim that $v(\cdot)$ is the convex conjugate of $u(\cdot)$, we define the auxiliary domain $\mathcal{A}'(x,\mathcal{E}) \triangleq \mathcal{A}(x,\mathcal{E}) \setminus \bigcup_{x'<x} \mathcal{A}(x',\mathcal{E})$. Note that (a) the monotonicity of the utility functional $\mathbf{U}(\cdot)$ implies that

$$\sup_{c\in\mathcal{A}(x,\mathcal{E})} \mathbf{U}(c) = \sup_{c\in\mathcal{A}'(x,\mathcal{E})} \mathbf{U}(c)$$

and (b) the Proposition 2.2 implies that $\sup_{\mathbb{Q}\in\mathcal{D}_\kappa(y)}\langle c-e, \mathbb{Q}\rangle = xy$ for any $y>0$ and $c \in \mathcal{A}'(x,\mathcal{E})$. Having established the weak-*compactness of the dual domain $\mathcal{D}_\kappa(y)$ in 3.4, the minimax theorem [see Sion (1958)] implies that

$$\begin{aligned}
\sup_{x\in\mathbb{R}}[u(x)-xy] &= \sup_{x\in\mathbb{R}}\bigg(\sup_{c\in\mathcal{A}'(x,\mathcal{E})} \mathbf{U}(c) - xy\bigg) \\
&= \sup_{x\in\mathbb{R}} \sup_{c\in\mathcal{A}'(x,\mathcal{E})} \bigg(\mathbf{U}(c) - \sup_{\mathbb{Q}\in\mathcal{D}_\kappa(y)}\langle c-e,\mathbb{Q}\rangle\bigg) \\
&= \sup_{x\in\mathbb{R}} \sup_{c\in\mathcal{A}'(x,\mathcal{E})} \inf_{\mathbb{Q}\in\mathcal{D}_\kappa(y)}(\mathbf{U}(c) - \langle c,\mathbb{Q}\rangle + \langle e,\mathbb{Q}\rangle) \\
&= \sup_{c\in\mathcal{V}_{\kappa+}^{\mathcal{M}}} \inf_{\mathbb{Q}\in\mathcal{D}_\kappa(y)}(\mathbf{U}(c) - \langle c,\mathbb{Q}\rangle + \langle e,\mathbb{Q}\rangle) \\
&= \inf_{\mathbb{Q}\in\mathcal{D}_\kappa(y)} \sup_{c\in\mathcal{V}_{\kappa+}^{\mathcal{M}}}(\mathbf{U}(c) - \langle c,\mathbb{Q}\rangle + \langle e,\mathbb{Q}\rangle) \\
&= \inf_{\mathbb{Q}\in\mathcal{D}_\kappa(y)}(\mathbf{V}(\mathbb{Q}) + \langle e,\mathbb{Q}\rangle) = v(y). \qquad \square
\end{aligned}$$

LEMMA A.2 (Existence in the dual problem). *For $y \in \mathrm{Dom}(v)$ there exists $\widehat{\mathbb{Q}}^y \in \mathcal{D}_\kappa(y)$ such that*

$$v(y) = \mathbf{V}^{\mathcal{E}}(\widehat{\mathbb{Q}}^y) = \mathbf{V}(\widehat{\mathbb{Q}}^y) + \langle e, \widehat{\mathbb{Q}}^y\rangle.$$

PROOF. For $y \in \mathrm{Dom}(v)$, let $(\mathbb{Q}_n)_{n\in\mathbb{N}}$ be a minimizing sequence for $v(y)$, that is, a sequence in $\mathcal{D}_\kappa(y)$, such that $(\mathbf{V}^{\mathcal{E}}(\mathbb{Q}_n))_{n\in\mathbb{N}}$ is real valued and decreasing with limit $v(y)$. Since $\mathcal{D}_\kappa(y)$ is a closed and bounded subset of the dual $(\mathcal{V}_\kappa^{\mathcal{M}})^*$ of $\mathcal{V}_\kappa^{\mathcal{M}}$, by Proposition 3.4, the product space $\mathcal{D}_\kappa(y) \times [v(y), \mathbf{V}^{\mathcal{E}}(\mathbb{Q}_1)]$ is compact. Therefore, the sequence $(\mathbb{Q}_n, \mathbf{V}^{\mathcal{E}}(\mathbb{Q}_n))_{n\in\mathbb{N}}$ has a cluster point $(\widehat{\mathbb{Q}}^y, v^*)$ in $\mathcal{D}_\kappa(y) \times [v(y), \mathbf{V}^{\mathcal{E}}(\mathbb{Q}_1)]$. By the decrease of the sequence $(\mathbf{V}^{\mathcal{E}}(\mathbb{Q}_n))_{n\in\mathbb{N}}$, we have $v^* = \lim_n \mathbf{V}^{\mathcal{E}}(\mathbb{Q}_n) = v(y)$. On the other hand, by the definition (4.1) of the functional $\mathbf{V}(\cdot)$, the epigraph of its restriction $\mathbf{V}^{\mathcal{E}}(\cdot) : \mathcal{D}_\kappa(y) \to \mathbb{R}$ is closed with respect to the product of the weak-* and Euclidean topologies. Therefore, $(\widehat{\mathbb{Q}}^y, v^*)$ is in the epigraph of $\mathbf{V}^{\mathcal{E}}$ and thus $v(y) = v^* \geq \mathbf{V}^{\mathcal{E}}(\widehat{\mathbb{Q}}^y) = \mathbf{V}(\widehat{\mathbb{Q}}^y) + \langle\widehat{\mathbb{Q}}^y, e\rangle$.  $\square$

LEMMA A.3 (Consequences of reasonable elasticity).

1. *We have* $\text{Dom}(v) = (0, \infty)$.
2. *We have $v(\cdot)$ is continuously differentiable and, for $y > 0$, its derivative satisfies*

$$yv'(y) = -\langle (\widehat{\mathbb{Q}}^y)^r, \mathbf{I}(\widehat{\mathbb{Q}}^y)\rangle + \langle e, \widehat{\mathbb{Q}}^y\rangle,$$

   *where $\widehat{\mathbb{Q}}^y \in \mathcal{D}_\kappa(y)$ is a minimizer in the dual problem [i.e., $v(y) = \mathbf{V}^{\mathcal{E}}(\widehat{\mathbb{Q}}^y)$].*
3. *The inequality*

$$yv'(y) \geq -\langle \mathbb{Q}^r, \mathbf{I}(\widehat{\mathbb{Q}}^y)\rangle + \langle e, \widehat{\mathbb{Q}}^y\rangle$$

   *holds for all $\mathbb{Q} \in \mathcal{D}_\kappa(y)$.*
4. *We have $\lim_{y \to 0} v'(y) = -\infty$ and $\lim_{y \to \infty} v'(y) \in [\inf_{\mathbb{Q} \in \mathcal{M}} \mathbb{E}^{\mathbb{Q}}[\mathcal{E}_T], \sup_{\mathbb{Q} \in \mathcal{M}} \mathbb{E}^{\mathbb{Q}}[\mathcal{E}_T]]$*
5. *We have $\mathbf{I}(\widehat{\mathbb{Q}}^y) \in \mathcal{A}(-v'(y), e)$ and $\langle \mathbf{I}(\widehat{\mathbb{Q}}^y), (\widehat{\mathbb{Q}}^y)^r\rangle = \langle \mathbf{I}(\widehat{\mathbb{Q}}^y), \widehat{\mathbb{Q}}^y\rangle$.*

PROOF. Thanks to the representation $v(y) = \mathbb{E} \int_0^T V(t, Y_t^{\widehat{\mathbb{Q}}^y}) \, d\kappa_t$, and the fact that $\mathbb{E} \int_0^T Y_t^{\mathbb{Q}} \, d\kappa_t \leq 1$ for all $\mathbb{Q} \in \mathcal{D}_\kappa(1)$, the proofs of parts 1–4 of this lemma follow (almost verbatim) the proofs of the following statements in Karatzas and Žitković (2003): 1. Lemma A.5, page 30; 2. Lemma A.6, page 31; 3. Proposition A.7, page 32. 4. Lemma A.8, page 33.

To prove claim 5, we observe that the combination of parts 3 and 4 implies that

$$\langle \mathbf{I}(\widehat{\mathbb{Q}}^y), y\mathbb{Q}\rangle \leq -yv'(y) + \langle e, y\mathbb{Q}\rangle \qquad \text{for all } \mathbb{Q} \in \mathcal{M}_\kappa.$$

From Proposition 3.5 it follows that $\mathbf{I}(\widehat{\mathbb{Q}}^y) \in \mathcal{A}(-v'(y), e)$, so $\langle \mathbf{I}(\widehat{\mathbb{Q}}^y), \mathbb{Q}\rangle \leq -yv'(y) + \langle e, \mathbb{Q}\rangle$ for all $\mathbb{Q} \in \mathcal{D}(y)$. In particular, $\langle \mathbf{I}(\widehat{\mathbb{Q}}^y), \widehat{\mathbb{Q}}^y\rangle \leq -yv'(y) + \langle e, \widehat{\mathbb{Q}}^y\rangle$, yielding immediately the inequality $\langle \mathbf{I}(\widehat{\mathbb{Q}}^y), \widehat{\mathbb{Q}}^y\rangle \leq \langle \mathbf{I}(\widehat{\mathbb{Q}}^y), (\widehat{\mathbb{Q}}^y)^r\rangle$. The second part of the claim follows by the trivial inequality $\langle \mathbf{I}(\widehat{\mathbb{Q}}^y), \widehat{\mathbb{Q}}^y\rangle \geq \langle \mathbf{I}(\widehat{\mathbb{Q}}^y), (\widehat{\mathbb{Q}}^y)^r\rangle$. □

LEMMA A.4 (Existence in the primal problem). *For $x > -\lim_{y \to \infty} v'(y)$, the primal problem (2.5) has a solution, that is, there exists $\hat{c}^x \in \mathcal{A}(x, \mathcal{E})$ such that $u(x) = \mathbf{U}(\hat{c}^x)$. Moreover, the optimal consumption density process $\hat{c}^x$ is $\mathbb{P}_\kappa$-a.s. unique.*

PROOF. Using the continuous differentiability of the dual value function $v(\cdot)$ and Lemma A.5, we conclude that for any $x > \lim_{y \to \infty} v'(y)$ there exists a unique $y > 0$ such that $v'(y) = -x$. Let $\widehat{\mathbb{Q}}^y$ be the solution to the dual problem that corresponds to $y$ and define the candidate solution $\hat{c}^x$ to the primal problem by

$$\hat{c}^x \triangleq \mathbf{I}(\widehat{\mathbb{Q}}^y).$$

By Lemma A.3, $\hat{c}^x \in \mathcal{A}(x, \mathcal{E})$. The optimality of the consumption density process $\hat{c}^x$ follows from the fact that

$$\mathbf{U}(\hat{c}^x) = \mathbf{U}(\mathbf{I}(\widehat{\mathbb{Q}}^y)) = \mathbf{V}(\widehat{\mathbb{Q}}^y) + \langle \mathbf{I}(\widehat{\mathbb{Q}}^y), \widehat{\mathbb{Q}}^y \rangle = \mathbf{V}(\widehat{\mathbb{Q}}^y) + \langle \mathbf{I}(\widehat{\mathbb{Q}}^y), (\widehat{\mathbb{Q}}^y)^r \rangle$$
$$= v(y) - yv'(y) = u(x),$$

using Lemma A.3 and the conjugacy of $u(\cdot)$ and $v(\cdot)$. The $\mathbb{P}_\kappa$-a.s. uniqueness of $\hat{c}^x$ is a direct consequence of the strict concavity of the mapping $x \mapsto U(\omega, t, x)$ coupled with convexity of the feasible set $\mathcal{A}(x, \mathcal{E})$. $\square$

LEMMA A.5. *We have* $\lim_{y \to \infty} v'(y) = \mathcal{L}(\mathcal{E})$, *where* $\mathcal{L}(\mathcal{E}) = \inf_{\mathbb{Q} \in \mathcal{M}} \mathbb{E}^\mathbb{Q}[\mathcal{E}_T]$.

PROOF. Let $x' = \lim_{y \to \infty} v'(y)$. Part 4 of Lemma A.3 states that $x' \geq \mathcal{L}(\mathcal{E})$, so we need only to prove that $x' \leq \mathcal{L}(\mathcal{E})$. Suppose, to the contrary, that there exists $x_0 > \mathcal{L}(\mathcal{E}_T)$ of the form $x_0 = v'(y_0)$ for some $y_0 > 0$ so that $x' > x_0$. The optimal consumption process $(C_t^{-x_0})_{t \in [0,T]}$ that corresponds to the initial capital $-x_0$ exists by Lemma A.4 and satisfies $\mathbb{E}^\mathbb{Q}[C_T^{-x_0}] \leq -x_0 + \mathbb{E}^\mathbb{Q}[\mathcal{E}_T]$ for any $\mathbb{Q} \in \mathcal{M}$ by Proposition 2.2. Taking the infimum over $\mathbb{Q} \in \mathcal{M}$, we reach a contradiction:

$$0 \leq \inf_{\mathbb{Q} \in \mathcal{M}} \mathbb{E}^\mathbb{Q}[C_T^{-x_0}] \leq -x_0 + \mathcal{L}(\mathcal{E}_T) < 0.$$

Therefore, $x' \leq \mathcal{L}(\mathcal{E})$. $\square$

**Acknowledgments.** I thank the anonymous referees for a number of useful suggestions and improvements. I am also indebted to the probability seminar participants at Carnegie Mellon University, Brown University and Boston University. Any opinions, findings, and conclusions or recommendations expressed in this material are my own and do not necessarily reflect the views of the National Science Foundation.

DEPARTMENT OF MATHEMATICAL SCIENCES
CARNEGIE MELLON UNIVERSITY
7209 WEAN HALL
PITTSBURGH, PENNSYLVANIA 15217
USA
E-MAIL: zitkovic@cmu.edu
URL: www.andrew.cmu.edu/~gordanz/